%% file: lbbdpaper.tex
\documentclass[review, authoryear]{elsarticle}

\usepackage[table]{xcolor}
\definecolor{gray1}{gray}{0.95}
\definecolor{gray2}{gray}{0.975}
\usepackage{lscape}

\usepackage{mathtools, amssymb, amsthm}

\usepackage[colorlinks=true, citecolor=red, breaklinks=true]{hyperref}

\journal{arXiv}
\newcommand{\ts}{\textsuperscript}
\newcommand{\st}{\text{subject to}}
\renewcommand{\leq}{\leqslant}
\renewcommand{\geq}{\geqslant}

\newcommand{\len}{\textsf{len}}
\newcommand{\Delay}{\textsf{Delay}}
\newcommand{\Mins}{\textsf{M}}
\newcommand{\CanStart}{\textsf{CanStart}}
\newcommand{\Avail}{\textsf{Avail}}

\usepackage{algorithm}
\usepackage[noend]{algpseudocode}

\newtheorem{lemma}{Lemma}[section]

\theoremstyle{definition}
\newtheorem{example}{Example}[section]

\usepackage[margin=1in]{geometry}

\usepackage{tikz}
\usepackage{pgfplots}
\usepackage{graphicx}
\graphicspath{ {./images/} }
\tikzset{>=stealth}

\usepackage{subcaption}
\DeclareCaptionFormat{custom}
{\textbf{#1#2}\textit{\small #3}}

\begin{document}

\begin{frontmatter}
\input{frontmatter.tex}

\end{frontmatter}

\noindent Accompanying GitHub:
\href{https://git.io/JWfBO}{https://git.io/JWfBO}

\noindent $^\dagger$Corresponding Author email: \href{mitchopt@gmail.com}{mitchopt@gmail.com}

\allowdisplaybreaks

\section{Introduction}\label{sec:1}\input{sec1.tex}

\section{Logic-Based Benders Decomposition}\label{sec:2}\input{sec2.tex}

\section{Benders Cuts for the Main Problem}\label{sec:3}\input{sec3.tex}

\section{A Class of Scheduling Problems}\label{sec:4}\input{sec4.tex}

\section{Nursing Home Shift Scheduling}\label{sec:5}\input{sec5.tex}

\section{Airport Check-in Counter Allocation}\label{sec:6}\input{sec6.tex}

\section{Conclusion}\label{sec:7}\input{sec7.tex}

\section*{Acknowledgements}\input{Acknowledgements.tex}

\addcontentsline{toc}{section}{References}
\bibliography{CombOptSimLBBD}

\appendix
\section{Implementation Details}\label{sec:8}\input{sec8.tex}

\end{document}

%% file: frontmatter.tex
\title{Combining Optimisation and Simulation Using Logic Based Benders Decomposition}

\author{M. A. Forbes}
\author{\quad M. G. Harris$^\dagger$}
\address{School of Mathematics and Physics, The University of Queensland, Brisbane, Australia}
\author{H. M. Jansen}
\address{Department of Engineering, University College Roosevelt, Middelburg, The Netherlands}
\author{F. A. van der Schoot}
\address{Department of Mathematics and Computer Science, Eindhoven University of Technology, Eindhoven, The Netherlands}
\author{T. Taimre}
\address{School of Mathematics and Physics, The University of Queensland, Brisbane, Australia}

\begin{abstract}
Operations Research practitioners often want to model complicated functions that are
are difficult to encode in their underlying optimisation framework. 
A common approach is to solve an approximate model, and then use a simulation to evaluate the true objective value of one or more solutions. 
We propose a new approach to integrating simulation into the optimisation model itself.
The idea is to run the simulation at each incumbent solution to a master problem.
The simulation data is then used to guide the trajectory of the optimisation model itself using logic based Benders cuts.
We test the approach on a class of stochastic resource allocation problems with monotonic performance measures.
We derive strong novel Benders cuts that are provably valid for all problems of the given form.
We consider two concrete examples: a nursing home shift scheduling problem, and an airport check in counter allocation problem.
While previous papers on these applications could only approximately solve realistic instances, we are able to solve them exactly within a reasonable amount of time.
Moreover, while those papers account for the inherent variance of the problem by including estimates of the underlying random variables as model parameters, we are able to compute sample average approximations to optimality with up to $100$ scenarios.
\end{abstract}

\begin{keyword}
	integer programming; logic-based Benders decomposition; simulation; resource allocation;
	shift scheduling
\end{keyword}

%% file: sec1.tex
In resource allocation problems, the goal is to optimally distribute scarce resources among a finite set of objects, such as time periods, or locations in a network.
Although the nature of the resources might vary from problem to problem, we are usually interested in how many resources should be allocated to each object, which makes integer programming a natural solution approach.
If the objective contribution of each object is a convex function of one variable, then standard techniques may be suitable. But in practice, objective functions are often non-convex, and include complex interactions between the variables.

The problems we study have the following mathematical form. Let $\mathbb{Z}_+$ and $\mathbb{R}_+$ denote the sets of non-negative integers and reals respectively. For a positive integer $n$, denote by $[n] = \{1, \ldots, n\}$ the set of integers ranging from $1$ up to $n$.
We study problems of the form
\begin{equation}
	\min_y \quad g(y) + \sum_{\omega \in \Omega}f_\omega(y) \quad \st \quad y \in Y, \label{problem:main}
\end{equation}
where $Y \subseteq \mathbb{Z}_+^n$ is a finite set, $g : \mathbb{Z}_+^n \to \mathbb{R}$ is function, and for $\omega \in \Omega$, $f_\omega : \mathbb{Z}_+^n \to \mathbb{R}_+$ is a bounded non-increasing function of each of the variables. The elements of $[n]$ represent a collection of objects, and for each $j \in [n]$, the $j\ts{th}$ component of $y = (y_1,\ldots,y_n)$ represents the number of resources allocated to the $j\ts{th}$ object. The objective function consists of a cost function, $g$, for the resources, and for each $\omega \in \Omega$, a performance measure, $f_\omega$. Our task is to find a resource allocation, $y$, which minimises the combined cost and performance. Our soft assumption is that $f_\omega$ are the outputs of a complicated simulation that depends on the $y$ variables. The importance of the monotonicity property for each $f_\omega$ will be made clear soon, but first, consider an example.

\begin{example}
	We are given a finite set of jobs, and $n$ adjacent time periods of finite length. Each job has a release time and a processing time. We can decide how many staff will work in each time period, subject to constraints, and the goal is to minimise the total delay (start time minus release time) of all of the jobs. The jobs are processed as staff become available, and are scheduled on a first-come-first-serve (FCFS) basis. For all $j \in [n]$, let $y_j$ be an integer variable representing the number of staff working in time period $j$, and let $\textsf{\textup{Delay}}_{j}(y)$ be the total delay of all jobs released in time period $j$, given the staffing vector $y = (y_1,\ldots,y_n)$. In Section \ref{sec:4} we will prove that $\textsf{\textup{Delay}}_{j}$ is a non-increasing function of all of the variables. In other words, increasing the number of staff available in \emph{any} time period cannot increase the delay of \emph{any} job. \hfill $\diamondsuit$ \label{example:intro}
\end{example}

This paper is motivated by several properties of this example. Firstly, suppose that at some point in time, not enough staff are available to start processing all of the jobs which have been released. Then we get a backlog of delayed jobs. So for each $j$, $\textsf{\textup{Delay}}_{j}(y)$ depends not only on the number of staff available in time period $j$ but on other time periods as well. It depends on previous time periods if the staff available in time period $j$ are busy processing the backlog, and on future time periods if the jobs released now are forced to join the backlog. Since the exact dependencies are not known a priori, $\textsf{\textup{Delay}}_{j}$ really is a function of the decision variables for each time period; not just $y_j$. Be that as it may, in practice the delays tend to depend strongly on neighbouring time periods, and weakly or not at all on distant ones. We will make this notion more precise later on.

Secondly, $\textsf{\textup{Delay}}_{j}$ does not have a closed-form analytic expression. It is possible to model the full problem as a single integer program, but incorporating the queuing structure of the delays requires a large number of new variables and constraints, and the resulting model is (as we will see) prohibitively difficult to solve. On the other hand, for a fixed choice of the staffing vector $y$, we can evaluate each $\textsf{\textup{Delay}}_{j}(y)$ during a discrete-event simulation (DES) of the job queue. A common approach is to solve an easier problem which approximates the queuing structure, and to use simulation modelling to evaluate the actual delay associated with one or more solutions. Simulating the performance of every feasible solution is clearly impractical, but thanks to the monotonicity property, the output of one simulation provides useful lower bounds on the performance of other solutions. These observations motivate us to consider Benders decomposition.

\subsection{Stochastic Optimisation}

When modelling real phenomena, the exact value of the performance measures frequently depends on information that cannot be known ahead of time. For instance, with Example \ref{example:intro}, we might be required to allocate the staff to their shifts before the exact release times, processing times, or even the number of jobs is known. To model this situation mathematically we assume $f_\omega$ also depends on one-or-more random variables, $\mathfrak R$. The problem to be solved is then
\begin{equation}
		\min_y \quad g(y) + \mathbb{E}\sum_{\omega \in \Omega}f_\omega(\mathfrak R, y) \quad \st \quad y \in Y, 
	\label{problem:stochastic}
\end{equation}
where $\mathbb{E}$ denotes the expected value operator. This is a stochastic integer programming problem, and trying to calculate the exact expectation is hopeless, because the performance measures lack any convenient analytic structure. A standard approach in this situation is to calculate the average performance over a finite number of samples. We either assume that $\mathfrak R$ can be sampled efficiently, or, if the context allows it, that we have access to historical data. To obtain a reliable approximation, we would like to use a relatively large number of samples, which makes it even more challenging to use a direct approach.

Let $S$ be the set of scenarios; that is, an indexing set for the samples. For each $s \in S$, let $r_{s} \sim \mathfrak R$ denote the $s\ts{th}$ sample, and define $f_{s\omega}(y) = f_\omega(r_s, y)$. The uniform sample-average approximation (SAA) of (\ref{problem:stochastic}) is then
\begin{equation}
		\min_y \quad g(y) + \sum_{s \in S}\sum_{\omega\in\Omega}f_{s\omega}(y)/\lvert S\rvert \quad \st \quad y \in Y, \label{problem:general}
\end{equation}
which is a deterministic optimisation problem. We assume $f_{s\omega}(y)$ is a non-increasing function in every scenario. Observe now that (\ref{problem:general}) is of the form (\ref{problem:main}) where $\Omega$ is replaced by $S \times \Omega$. Note also that (\ref{problem:general}) no longer depends explicitly on $\mathfrak R$. Therefore the approach is independent of certain distributional assumptions that may limit the applicability of the model. The only requirement is that we can (at least approximately) sample from the underlying distributions. Later it will be useful to assume that each of the performance measures is associated with a unique object, which we express using a function $\pi : \Omega \to [n]$. In other words, $\pi(\omega)$ is the object associated with the $\omega\ts{th}$ performance measure. For instance, in Example \ref{example:intro}, $\Omega = [n]$, and the object most naturally associated with $\textsf{Delay}_j$ is $j$ itself.

\subsection{Literature Review}

Classical Benders decomposition was introduced by \cite{Benders1962}, and we refer the reader to \cite{Rahmaniani2017} for a survey of applications. \cite{Hooker2000, HookerOttosson2003} introduced \textit{logic-based Benders decomposition} (LBBD), an ambitious generalisation of classical Benders decomposition. They observed that the concept underlying Benders decomposition could be applied in much more general settings. Since then, LBBD has been applied with remarkable success to a diverse range of problems. We refer the reader to \cite{Hooker2019} for a survey of applications.
Like classical Benders decomposition, the original problem is partitioned into a master problem and a set of subproblems. The master problem is a relaxation of the original problem, and we use new variables to approximate the original objective function. At incumbent solutions to the master problem, the subproblems are solved to generate Benders cuts which gradually refine the approximation. In our case, the resource allocation decisions will be left in the master problem, and the computation of the performance measures will be relegated to the subproblems. The master problem is a mixed integer program (MIP), and we solve the subproblems using discrete-event simulations.

To apply classical Benders decomposition, the subproblem must be a linear program (LP), and cuts are derived by solving the dual problem. \cite{Geoffrion1972} introduced Generalized Benders decomposition, which is an extension to certain non-linear programming problems that uses convex duality theory instead. Another variant, Combinatorial Benders decomposition, was introduced by \cite{CodatoFischetti2006} and is suited to problems with implication constraints. With LBBD on the other hand, the subproblems can take the form of any optimisation problem. In fact, the subproblem can be any function evaluation, as long as suitable Benders cuts can be derived.
While this flexibility has allowed the principles of Benders decomposition to be applied to a larger set of problems, one weakness is that the structure of the Benders cuts is problem-specific. Indeed, finding suitable Benders cuts sometimes requires considerable effort, since we lack a convenient theory of duality for an arbitrary subproblem. In the case of (\ref{problem:main}), our assumptions are weak enough to capture a large number of practical problems. Nevertheless, we can derive strong Benders cuts which are provably valid for the entire class of problems.

Classical Benders decomposition has been applied to two-stage stochastic programming problems with linear recourse by \cite{SlykeWets1969}, with integer recourse by \cite{LaporteLouveaux1993}, and in many other papers following those studies. In the case of integer recourse, Benders cuts are derived by solving the dual of the LP relaxation. To ensure convergence they are supplemented with combinatorial cuts that eliminate previous master solutions. This approach suffers from slow convergence if the subproblem has a weak linear relaxation, since combinatorial cuts may be needed at many master solutions. Logic-based Benders decomposition, on the other hand, does not rely on LP duality to derive Benders cuts. \cite{OzgunHooker2020} apply LBBD to two-stage stochastic planning and scheduling problems, where the master problem is an assignment problem that is solved using mixed-integer programming, and the subproblems are scheduling problems, on which constraint programming excels. Benders cuts are derived by solving the so-called ``inference dual'' of the subproblems. We refer the reader to \cite{HookerOttosson2003} for an explanation of the inference dual.

\cite{LombardiMilanoRuggieroBenini2010} use LBBD to solve a \textit{Stochastic Allocation and Scheduling} problem for conditional task graphs in multi-processor embedded systems. The master problem assigns tasks to chips, and tasks are scheduled by the subproblems. The subproblems communicate with the master problem via so-called ``no-good cuts,'' which we will explain in Section \ref{sec:2}. \cite{ZarandiBermanBeck2012} apply LBBD to a \textit{Stochastic Facility Location and Vehicle Assignment} problem, where the subproblems communicate with the master problem via feasibility cuts. \cite{GuoBodurAlemanUrbach2019} apply LBBD to a stochastic extension of the \textit{Distributed Operating Room Scheduling} problem, which aims to find an assignment of surgeries to operating theatres. In that problem, the subproblems are always feasible, and optimality cuts are used instead. We will explain the difference between feasibility and optimality cuts in Section \ref{sec:2}.

Monotonicity is a strong property which has contributed to logical cuts elsewhere in the literature. For instance, in the context of job shop scheduling, adding a job to a machine can only increase the makespan for that machine (and removing a job can only reduce it); see \cite{NaderiRoshanaei2021} for instance. In the context of stochastic operating room scheduling, \cite{GuoBodurAlemanUrbach2019} exploits the monotonicity of the surgery cancellation costs. In the context of facility location and vehicle assignment, \cite{ZarandiBermanBeck2012} exploit the fact that increasing the number of clients assigned to a facility cannot reduce the number of trucks required. \cite{Fischetti2019} use a ``downward monotonicity'' property to generate cuts for certain two-person interdiction games. There the assumption is that reducing the number of interdictions preserves the feasibility of a follower policy. In contrast, we exploit a stronger form of monotonicity in multiple variables; recall Example \ref{example:intro}, where we note that increasing the number of staff at \textit{any} time cannot increase the delay of \textit{any} job.
  
It is well known that classical Benders decomposition excels if the subproblems can be solved analytically or using an efficient special-purpose algorithm, rather than with a generic LP solver. \cite{PearceForbes2018} apply classical Benders decomposition to a \textit{Dynamic Facility Location and Network Design} problem. They show that the dual variables associated with the subproblems can be calculated analytically, and, therefore, that incumbent solutions to the master problem can be evaluated almost instantly. \cite{ZhangMattaAlfieriPedrielli2017} apply classical Benders decomposition to a \textit{Joint Workstation, Workload, and Buffer Allocation} problem, which deals simultaneously with the design of an open flow line, and the total processing time to be allocated to each station. It is observed that the optimal dual variables of the linear subproblem can be calculated directly by simulating the performance of the system given the optimal master solution. The same group applies a similar approach in \cite{ZhangMattaAlfieriPedrielli2018, ZhangMattaAlfieriPedrielli2019}.
In this paper, we show that simulation can also be used to derive logic-based Benders cuts, where no dual variables are available.

In Section \ref{sec:5} we will study our first application; a \textit{Nursing Home Shift Scheduling} (NHSS) problem. In Section \ref{sec:6} we study an \textit{Airport Check-in Counter Allocation} (ACCA) problem. Application-specific literature reviews for these problems have been deferred to the relevant sections. Most existing approaches in these areas have the following features: the stochastic nature of the problem is accounted for by including estimates of the underlying variables as model parameters, and the queuing structure of the problem is only approximately modelled. Benders decomposition has a number of advantages. Since the simulation data will inform the optimization model, no approximation of the queuing structure is required. At the level of the scenarios, we get an exact solution. And because we can solve to optimality over a large number of scenarios, we can be confident that the solution closely approximates the true optimal solution to the stochastic problem.

\subsection{Outline of the Paper}

In Section \ref{sec:2} we provide a mathematical outline of the LBBD method. In Section \ref{sec:3} we derive a sequence of Benders cuts which are valid for all problems that can be written in the form (\ref{problem:main}), and two families of valid inequalities which can be added to the master problem in advance. In Section \ref{sec:4} we expand on Example \ref{example:intro} in more detail, and prove that it has the form of (\ref{problem:main}). In Section \ref{sec:5} we will study  the Nursing Home Shift Scheduling (NHSS) problem. In Section \ref{sec:6} we study the Airport Check-in Counter Allocation (ACCA) problem. In Section \ref{sec:7} we conclude with a discussion of the results and some directions for future research.

%% file: sec2.tex
To set the stage for the problem we are interested in, in this section we will provide a general mathematical framework for the LBBD method.
Consider an optimisation problem of the following form, which we call the \emph{original problem} (OP):
\begin{equation}
		\min_y \quad g(y) + \sum_{\omega \in \Omega}\Theta_{\omega}(y) \quad \st \quad y \in Y,
\end{equation}
where $Y, Y'$, and $\Omega$ are sets,
and $g : Y' \to \mathbb{R}\cup\{\pm\infty\}$ and $\Theta_{\omega} : Y' \to \mathbb R \cup\{\pm\infty\}$ for $\omega \in \Omega$ are extended real valued functions.
Furthermore, $Y$ is the feasible region, $Y'$ is the domain of the objective function, and $\Omega$ is finite.
We adopt the convention whereby a minimisation problem has $+\infty$ objective value if it is infeasible, and $-\infty$ if it is unbounded.
Intuitively speaking, we think of $g$ as the \textit{easy} part of the objective function, while $\Theta_\omega$ are the \textit{hard} parts, such as:
\begin{enumerate}
	\item The value function of another minimisation problem (the \emph{subproblem}) parameterised by the $y$ variables.
	\item An indicator function for a logical proposition; in other words, $\Theta_\omega(y) = \infty$ if $P(y)$ is true, and $\Theta_\omega(y) = 0$ if $P(y)$ is false.
	This enforces that the logical proposition $P$ is true of any feasible solution.
	\item The output of a measurement, simulation, or function evaluation.
\end{enumerate}

In other words, we assume that the problem reduces to a significantly easier one if we ignore the $\Theta_\omega$ functions.
For simplicity we will refer to $\Theta_{\omega}$ as subproblems whether or not they correspond to optimisation problems per se.
An important special case occurs if $Y$ is the intersection of a polyhedron with the integer lattice, and $g$ is a linear function. In that case, the problem to be solved is \emph{almost} an integer program (IP) save for the inclusion of the complicating $\Theta_\omega$ functions. For example, in classical Benders decomposition, $\Theta_\omega$ is the value function of an LP parameterised by the integer variables. Traditionally the integer variables have been thought of as the complicating factors. But although LPs are easier to solve than IPs in both theory and practice, it is that we carry an LP \textit{for each} $y \in Y$ that complicates what is otherwise a tractable integer program.

To apply Benders decomposition to OP we first ask ourselves the following question:
\begin{itemize}
	\item[] For each $\omega \in \Omega$ and $y' \in Y$, can we find a function $B_{\omega y'} : Y \to \mathbb{R}$  such that
	\begin{equation}
	B_{\omega y'}(y) \leq \Theta_{\omega}(y)\label{eq:benders_qn}
	\end{equation}
	is a valid inequality that obtains equality at $y = y'$, where $B_{\omega y'}$ can be encoded in our underlying optimisation framework?
\end{itemize}
If so, then the problem is a candidate for Benders decomposition. If $\Theta_\omega(y') = \infty$, then by convention, the subproblem associated with $\omega$ and $y'$ is infeasible. Therefore $y'$ cannot be an optimal solution to OP, and, in that case, it is enough to find an inequality which strictly separates $y'$ from the set of feasible solutions. To apply Benders decomposition in practice, we first replace the $\Theta_\omega$ functions with a new collection of continuous variables, $\theta_\omega$, intended to estimate their contribution to the objective function. The resulting relaxation is called the \emph{master problem} (MP) and it has the following form:
\begin{equation}
\min_y \quad g(y) + \sum_{\omega \in \Omega}\theta_{\omega} \quad \st \quad \theta \in \mathbb{R}^{\Omega},\, y \in Y.
\tag{\textrm{MP}}
\end{equation}

If the range of $\Theta_\omega$ is $\{0, \infty\}$ then we can omit the corresponding variable. By assumption, MP is significantly easier to solve than OP and yields a lower bound on its optimal objective value. Currently this is obvious since MP is unbounded, but we make sure the property is preserved by Benders cuts. For now let $\theta^*,\, y^*$ be an optimal solution to MP in its current form. We solve the subproblem associated with this solution by evaluating $\Theta_{\omega}(y^*)$ for each $\omega \in \Omega$, whether that is solving an optimisation problem, verifying a logical proposition, running a simulation, or something else. Then we consider the following case distinctions:
\begin{enumerate}
	\item If $\theta_\omega^* = \Theta_{\omega}(y^*) < \infty$ for all $\omega \in \Omega$, then $y^*$ is an optimal solution to OP.
	\item For each $\omega \in \Omega$ such that $\theta_{\omega}^* < \Theta_{\omega}(y^*) < \infty$, add a Benders optimality cut of the form 
	\begin{equation}
	B_{\omega y^*}(y) \leq \theta_\omega
	\end{equation}to MP and continue. A trivial Benders cut is given by
	\begin{equation*}
		B_{\omega y^*}(y) = \begin{cases}
			\Theta_{\omega}(y^*) & \text{if }y = y^*, \\
			-\infty & \text{otherwise}.
		\end{cases}
	\end{equation*}
	\textit{Cuts of this form are called} no-good cuts \textit{since they correct the underestimate at the current ``no good'' solution, but} only \textit{at that solution.}
	\item For each $\omega \in \Omega$ such that $\Theta_\omega(y^*) = \infty$, find a function $F_{\omega y^*}$ such that $F_{\omega y^*}(y)\leq 0$ if $y$ is feasible for OP, but $F_{\omega y^*}(y^*) > 0$.
	Then add a Benders \textit{feasibility} cut of the form 
	\begin{equation*}
	F_{\omega y^*}(y) \leq 0
	\end{equation*}
	to MP and continue. \textit{A no-good feasibility cut separates} $y^*$ \textit{from} $Y\setminus\{y^*\}$. \textit{In other words, it eliminates the current solution, but} only \textit{that solution}.
\end{enumerate}

Benders Decomposition iterates between solving the master problem and evaluating solutions. When the $\theta_\omega$ variables correctly estimate the subproblem values at an optimal master solution, we can stop. Otherwise we add optimality and feasibility cuts as necessary, and solve the new master problem. We continue in this way until the subproblem solutions agree with the master estimates.

\cite{Erlendur2001} introduced \textit{Branch-and-Check}, which is a variant of LBBD that evaluates each incumbent solution to a single master problem; in other words, Benders cuts are embedded into a Branch-and-Cut framework. While Branch-and-Check has not been shown to universally outperform standard LBBD (see \cite{Beck2010}) it often enjoys a number of advantages. Since rebuilding a Branch-and-Bound tree from scratch in each iteration is often cumbersome and wasteful, in practice, we prefer to add Benders cuts as lazy constraints during a callback routine. On the other hand, standard LBBD may be preferable when the master problem is much easier to solve than the subproblems. In our implementations, we use Branch-and-Check.

By definition, a Benders cut must be tight at the current solution, and valid at all others. In the worst-case scenario, we can use no-good cuts. Logic-based Benders decomposition gets its name from the important role of logical inference in deriving stronger cuts. The strength of an optimality cut is determined by how tight the inequality (\ref{eq:benders_qn}) is at other solutions, while the strength of a feasibility cut corresponds to how many other infeasible solutions it eliminates. The task of the practitioner is to \emph{infer} from the current subproblem solution, valid bounds on as many other solutions as possible.

When MP is a MIP, it is easy to see why stronger cuts translate into a faster Branch-and-Check algorithm. Consider the Branch-and-Bound tree of the master problem. A priori, the best known lower bound for each solution is trivial. Therefore, with no-good cuts, we have to build the entire tree in order to prove an optimal solution. But if our cuts impose non trivial bounds at other solutions as well, then we already have useful lower bounds at other nodes. Therefore, we may be able to fathom some nodes without explicitly enumerating all of their child nodes. The more solutions we can infer valid bounds for at a time, and the tighter those bounds are, the faster we can prune the Branch-and-Bound tree and prove optimality for MP. Similar reasoning also applies to a standard LBBD algorithm.

In the next section, we will derive valid optimality cuts for our class of problems.

%% file: sec3.tex
In this section we will derive valid Benders optimality cuts for any problem that can be written in the form of (\ref{problem:main}).
To do this, we first need to convert the integer variables into binary ones.
For each $j \in [n]$ and $\xi \in \{y_{\min},\ldots, y_{\max}\}$, introduce a binary variable $z_{j\xi}$ 
with $z_{j\xi} = 1$ if and only if $y_j = \xi$.
In other words, $z_{j\xi} = 1$ if $\xi$ resources are allocated to the $j\ts{th}$ object, and $z_{j\xi} = 0$ otherwise.
To connect the $y$ and $z$ variables we add the following linear constraints to MP:
\begin{equation}
	\sum_{\xi = y_{\min}}^{y_{\max}} z_{j\xi} = 1 \quad \text{and} \quad \sum_{\xi = y_{\min}}^{y_{\max}}\xi z_{j\xi} = y_j \quad \forall j \in [n]. \label{linkYtoZ}
\end{equation}
Let $\left(\theta', y', z'\right)$ be the incumbent integer solution to the master problem in its current form, with $\theta'$ optimal given $y'$. We solve the subproblem by evaluating $f_\omega(y')$ for all $\omega \in \Omega$. Suppose there exists $\omega \in \Omega$ with $\theta_{\omega}' < f_{\omega}(y')$. The
no-good cut on $(\omega, y')$ is
\smallskip
	\begin{equation} 
	\theta_\omega \geq f_\omega(y')\left(1 - \sum_{j \in [n]}\sum_{\xi \,=\, y_{\min}}^{y_j' - 1} z_{j\xi} \;-\; \sum_{j \in [n]}\sum_{\xi \,=\, y_j' + 1}^{y_{\max}} z_{j\xi}\right). \label{cut:1}
	\end{equation}

\smallskip
\noindent To see that (\ref{cut:1}) is a valid optimality cut, note that the expression in parentheses is equal to $1$ if and only if $y = y'$. Otherwise, the expression is non positive, and the inequality is trivial.

We can improve the no-good cut by exploiting the monotonicity of the performance measures. A priori, the only solution we know $f_\omega(y')$ to be a valid bound for is $y'$ itself, but since $f_\omega$ is non-increasing, if we decrease the number of resources allocated to any object, then $f_\omega$ cannot decrease. Therefore the monotonic cut,
	\begin{equation} 
	\theta_\omega \geq f_\omega(y')\left(1 - \sum_{j \in [n]}\sum_{\xi \,=\, y_j' + 1}^{y_{\max}} z_{j\xi}\right) \label{cut:2}
	\end{equation}
is also a valid optimality cut.

Recall that while $f_\omega$ is a function of all of the decision variables, in practice it might depend only on a proper subset of the objects. In Example \ref{example:intro}, for instance, the delay of the jobs released in a given time period depends mostly on the number of staff in nearby time periods. If there are not too many concurrent jobs, then the delays are unlikely to depend at all on distant time periods. This is contingent on the current solution; the fewer staff which allocated to adjacent time periods, the wider the window of dependency will be. To exploit this more generically we will need some new notation.

Given $y \in Y$ and any subset $\mathcal L \subseteq [n]$, we let $\Delta(\mathcal L, y) \in \mathbb{Z}_+^n$ denote the integral vector which we get by increasing $y_j$ to $y_{\max}$ for all $j \in [n]\setminus \mathcal L$.  For example,
%
\begin{equation*}
\Delta(\{3,4\},(4,3,3,2,4)) = (y_{\max},y_{\max},3,2,y_{\max}).
\end{equation*}
In other words, we specify the number of resources allocated to objects in $\mathcal L$, and make the most optimistic assumptions possible elsewhere. 
Since $f_\omega$ is non-increasing, we have $f_\omega(\Delta(\mathcal L, y)) \leq f_\omega(y)$ for all $y \in Y$ and all $\mathcal L \subseteq [n]$. 
If equality is obtained, however, then we have shown that $f_\omega$ does not depend on the number of resources allocated to objects in $[n]\setminus \mathcal L$, as long as we do not increase the numbers \textit{in} $\mathcal L$. 
In other words, if $f_\omega(\Delta(\mathcal L, y')) = f_\omega(y')$, then the local cut,
%
\begin{equation}
	\theta_\omega \geq f_\omega(y')\left(1 - \sum_{j \in \mathcal L}\sum_{\xi \,=\, y_j' + 1}^{y_{\max}} z_{j\xi}\right) \label{cut:3}
\end{equation}
is a valid optimality cut. This cut is stronger than the previous one if $\mathcal L$ is a proper subset of $[n]$, and it is much stronger if $\mathcal L$ is small.
Thus we have introduced a trade-off between the benefit of finding a small set $\mathcal L$ which obtains equality and the computational expense of doing so. As there is no way to derive a minimal set $\mathcal L$ analytically, we will find good candidates using a search algorithm. In practice, this means simulating the performance of a several intermediate solutions.

Cuts (\ref{cut:1}), (\ref{cut:2}), and (\ref{cut:3}) function by multiplying the desired bound by a logical expression which deactivates the cut if we cannot be sure it is valid at the current solution. 
To strengthen (\ref{cut:3}), we would also like to make the coefficients of the binary variables as small as possible. 
Indeed, if we run a few more simulations, we can bound the error in the $f_\omega(y')$ estimate when the number of resources \textit{are} increased.
We distinguish between increasing the number of resources allocated to $\pi(\omega)$, and to any other object in $\mathcal L$, but mathematically speaking we could ``centre'' the cut on any of the objects.

Let $e_j = (0,\ldots,1,\ldots,0) \in \mathbb{R}^{n}$ be the $j\ts{th}$ standard basis vector with a $1$ in the $j\ts{th}$ position, and $0$s elsewhere. 
For each $\xi \in \{y_{\min}, \ldots, y_{\max}\}$ we define an expression, $I_\omega(\xi)$, which is the performance of $\omega$ if we allocate $\xi$ resources to $\pi(\omega)$, and the maximum number of resources elsewhere. Formally,
\begin{equation*}
	I_\omega(\xi) = f_\omega(\Delta(\{\pi(\omega)\}, \xi e_{\pi(\omega)})).
\end{equation*}
Recall that $\pi(\omega)$ is an object naturally associated with the $\omega\ts{th}$ performance measure. Our assumption is that $f_\omega$ depends very strongly on the number of resources allocated to $\pi(\omega)$.
Since $f_\omega$ is non-increasing, if $\xi \geq y_{\pi(\omega)}'$, then $I_\omega(\xi)$ is a valid lower bound on $\theta_\omega$ if we increase the number of resources allocated to $\pi(\omega)$ from $y_{\pi(\omega)}'$ to $\xi$, no matter what we do elsewhere. 

Now define
\begin{equation*}
	\textsf{Base}_\omega(y') = f_\omega(\Delta(\{\pi(\omega)\}, y_{\pi(\omega)}'e_{\pi(\omega)})).
\end{equation*}
Again by monotonicity, $\textsf{Base}_\omega(y')$ is a valid lower bound on $\theta_\omega$ if we \emph{do not} increase the number of resources allocated to $\pi(\omega)$, no matter what we do elsewhere.
This proves that the strengthened cut,
%
%
\begin{equation}
		\theta_\omega \geq \; f_\omega(y') 
		\;- \; \sum_{\xi \,=\, y_{\pi(\omega)}' + 1}^{y_{\max}} (f_\omega(y') - I_\omega(\xi)) z_{\pi(\omega)\xi}
		- \; \sum_{j \in \mathcal L\setminus\{\pi(\omega)\}}
		\; \sum_{\xi \,=\, y_j' + 1}^{y_{\max}} \left(f_\omega(y') - \textsf{Base}_\omega(y')\right)z_{j\xi}
	\label{cut:4}
\end{equation}
%
%
is a valid Benders cut. Instead of deactivating the cut, we compensate the bound by an appropriate amount when the number of resources allocated to objects in $\mathcal L$ are increased. The new bounds might be weaker, but if they are not trivial, they may still save us explicitly enumerating some solutions. We are not prevented from increasing the resources both at $\pi(\omega)$ and and other elements of $\mathcal L$ at the same time, but then we suffer both penalties. If we increase the number of allocated resources for too many of the objects at once, the bound may become effectively trivial. Nevertheless, this cut is at least as strong as the previous three. The comparative strength of the four Benders cuts is illustrated in theory by the following example, and in practice by the computational results of Section \ref{sec:6}.

	\begin{example} \label{example:cut_strength} Let $n = 5$, $y_{\min} = 1$, and $y_{\max} = 5$, with no further constraints. Then there are $3125$ feasible solutions.
	Let $y' = (4, 3, 3, 2, 4)$ be the incumbent solution. 
	Cut (\ref{cut:2}) is non-trivial at $288$ solutions. 
	But if $\mathcal L = \{3,4\}$ and $f(y') = f(\Delta(\{3,4\}, y'))$, then (\ref{cut:3}) is non-trivial at $5 \times 5 \times 3 \times 2 \times 5 = 750$ different solutions, and tight at at least $5^3 = 125$ of those solutions. Cut (\ref{cut:4}) provides weaker but non-trivial bounds at additional solutions.
	\end{example}
\begin{figure}[h]
	\begin{center}
			\begin{tikzpicture}[scale=0.75]
					\foreach \x in {1,...,5} 
					\foreach \y in {1,...,5} {\fill[opacity = 0.75] (\x - 1, \y) circle (2pt);}
					\foreach \x in {1,...,5} {\node at (\x - 1, .5) {\x};}
					\node at (-1, 1) {$y_{\min}$};
					\node at (-1, 5) {$y_{\max}$};
					\draw[line width = 1.5pt] (-0.5, 4)--(0.5, 4)--(0.5, 3)
					--(2.5, 3)--(2.5, 2)--(3.5, 2)--(3.5, 4)--(4.5, 4);
					\fill (0, 4) circle (2.5pt);
					\fill (1, 3) circle (2.5pt);
					\fill (2, 3) circle (2.5pt);
					\fill (3, 2) circle (2.5pt);
					\fill (4, 4) circle (2.5pt);
				\end{tikzpicture} \quad
			\begin{tikzpicture}[scale=0.75]
					\foreach \x in {1,...,5} 
					\foreach \y in {1,...,5} {\fill[opacity = 0.75] (\x - 1, \y) circle (2pt);}
					\foreach \x in {1,...,5} {\node at (\x - 1, .5) {\x};}
					\draw[line width = 1.5pt] (-0.5, 4)--(0.5, 4)--(0.5, 3) 
					--(2.5, 3)--(2.5, 2)--(3.5, 2)--(3.5, 4)--(4.5, 4);
					\fill[gray, opacity = .25] (-0.5, 1) rectangle (0.5, 4);
					\fill[gray, opacity = .25] (0.5, 1) rectangle (2.5, 3);
					\fill[gray, opacity = .25] (2.5, 1) rectangle (3.5, 2);
					\fill[gray, opacity = .25] (3.5, 1) rectangle (4.5, 4);
					\fill (0, 4) circle (2.5pt);
					\fill (1, 3) circle (2.5pt);
					\fill (2, 3) circle (2.5pt);
					\fill (3, 2) circle (2.5pt);
					\fill (4, 4) circle (2.5pt);
				\end{tikzpicture} \quad
			\begin{tikzpicture}[scale=0.75]
					\foreach \x in {1,...,5} 
					\foreach \y in {1,...,5} {\fill[opacity = 0.75] (\x - 1, \y) circle (2pt);}
					\foreach \x in {1,...,5} {\node at (\x - 1, .5) {\x};}
					\draw[line width = 1.5pt] (-0.5, 4)--(0.5, 4)--(0.5, 3)
					--(2.5, 3)--(2.5, 2)--(3.5, 2)--(3.5, 4)--(4.5, 4);
					\fill[gray, opacity = .25] (-0.5, 1) rectangle (0.5, 5);
					\fill[gray, opacity = .25] (0.5, 1) rectangle (1.5, 5);
					\fill[gray, opacity = .25] (1.5, 1) rectangle (2.5, 3);
					\fill[gray, opacity = .25] (2.5, 1) rectangle (3.5, 2);
					\fill[gray, opacity = .25] (3.5, 1) rectangle (4.5, 5);
					\fill (0, 4) circle (2.5pt);
					\fill (1, 3) circle (2.5pt);
					\fill (2, 3) circle (2.5pt);
					\fill (3, 2) circle (2.5pt);
					\fill (4, 4) circle (2.5pt);
				\end{tikzpicture} 
			\caption{The no-good (left), monotonic (centre), and local cuts (right). The lattice points correspond to the binary variables. The solutions that the bound is valid at correspond to taxi-cab-like paths through the shaded region.} \label{fig:cutcompare}
		\end{center}
\end{figure}
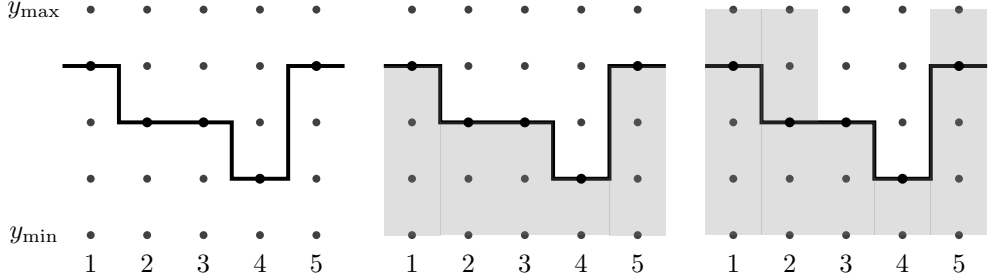
\vspace*{-2em}

\subsection{Initial Cuts}

We can get useful initial cuts by considering what happens if we specify number of resources allocated to one or two objects at a time, and make the most optimistic assumption possible elsewhere. In the language of the LBBD literature, we will add a relaxation of the subproblem to the master problem. First of all, we know that
\begin{equation}
	\theta_\omega \geq \sum_{\xi \,=\, y_{\min}}^{y_{\max}}I_\omega(\xi)z_{\pi(\omega)\xi} \quad \text{ for all } \omega\in\Omega
\end{equation}
is valid for the master problem. This is because $I_\omega(\xi)$ is a valid bound on $f_k$ whenever $\xi$ resources are allocated to $\pi(\omega)$, and precisely one term on the right-hand side is non-zero. This cut describes the convex hull of $I_\omega$ on $\{y_{\min}, \ldots, y_{\max}\}$, which makes it very valuable in the LP relaxation. The initial bounds can possibly be strengthened if we consider a second object at a time. Let $k \in [n]$ and $l \in [n]\setminus\{k\}$ be a pair of objects. For $\xi_1 \in \{y_{\min}, \ldots, y_{\max}\}$ and $\xi_2 \in \{y_{\min}, \ldots, y_{\max}\}$, let
\begin{equation}
	W_{kl}^\omega(\xi_1, \xi_2) = f_\omega(\Delta(\{k, l\}, \xi_1 e_{k} + \xi_2 e_{l}))
\end{equation}
denote the performance of $\omega$ whenever $\xi_1$ resources are allocated to $k$, and $\xi_2$ resources are allocated to $l$. Then
\begin{equation}
	\theta_\omega \geq W_{\pi(\omega)l}^\omega(\xi_1, \xi_2)(z_{\pi(\omega)\xi_1} + z_{l\xi_2} - 1)
\end{equation}
is valid for the master problem, having put $k = \pi(\omega)$. Computational experience suggests that this cut is weak, which makes sense, since the cut is only active if the correct number of resources are allocated to both objects at once. As an intermediate step towards strengthening it, fix $\xi_1$, and suppose we tried to impose the bound of $W_{\pi(\omega)l}^\omega(\xi_1, \xi_2)$ on $\theta_\omega$ whenever $\xi_2$ resources are allocated to $l$, regardless of how many are allocated to $\pi(\omega)$. The necessary inequality would be 
\begin{equation*}
	\theta_\omega \geq \sum_{\xi_2 = y_{\min}}^{y_{\max}}W_{\pi(\omega)l}^\omega(\xi_1, \xi_2)z_{l\xi_2}(y)
\end{equation*}
which is stronger, but no longer valid in general. However, since $f_\omega$ is non-increasing, it \emph{is} valid if $y_{\pi(\omega)} \leq \xi_1$. Suppose, on the other hand, that $y_{\pi(\omega)} > \xi_1$. Then regardless of the value of $\xi_2$, the error in our estimate is no greater than
\begin{equation*}
	\max_{\xi_2 \in \{y_{\min}, \ldots, y_{\max}\}}\left( W_{\pi(\omega)l}^\omega(\xi_1, \xi_2) - W_{\pi(\omega)l}^\omega(y_{\pi(k)}, \xi_2) \right).
\end{equation*}
It follows that
\begin{align}
	\theta_\omega \geq &\sum_{\xi_2 = y_{\min}}^{y_{\max}}W_{\pi(\omega)l}^\omega(\xi_1, \xi_2)z_{l\xi_2}(y) \nonumber \\
	& - \sum_{\xi' \,=\, \xi_1 + 1}^{y_{\max}} \max_{\xi_2 \in \{y_{\min},\ldots, y_{\max}\}}\left( W_{\pi(\omega)l}^\omega(\xi_1, \xi_2) - W_{\pi(\omega)l}^\omega(\xi', \xi_2) \right)z_{\pi(\omega)\xi'}(y)
	\label{initial:1}
\end{align}
is valid for the master problem for all $\xi_1 \in \{y_{\min}, \ldots, y_{\max}\}$.

An analogous argument, beginning instead with ``suppose we tried to impose the bound of $W_{\pi(\omega)l}^\omega(\xi_1, \xi_2)$ on $\theta_\omega$ whenever $\xi_1$ resources are allocated to $\pi(\omega)$, regardless of how many are allocated to $l$,'' leads us to a second set of cuts:
\begin{align}
	\theta_\omega \geq & \sum_{\xi_1 = y_{\min}}^{y_{\max}}W_{\pi(\omega)l}^\omega(\xi_1, \xi_2)z_{\pi(\omega)\xi_1}(y) \nonumber \\ & -
	\sum_{\xi' \,=\, \xi_2 + 1}^{y_{\max}}
	\max_{\xi_1 \in \{y_{\min}, \ldots, y_{\max}\}}\left( W_{\pi(\omega)l}^\omega(\xi_1, \xi_2) - W_{\pi(\omega)l}^\omega(\xi_1, \xi') \right)z_{l\xi'}(y)
	\label{initial:2}
\end{align}
for all $\xi_2 \in \{y_{\min}, \ldots, y_{\max}\}$. The one-dimensional cuts described earlier are the special cases which we obtain when either $\xi_1 = y_{\max}$ or $\xi_2 = y_{\max}$, so only the two-dimensional cuts will be added in practice.

A pair of cuts are formally valid for the master problem for all $\omega \in \Omega$ and $l \in [n]\setminus \pi(\omega)$.  But generating all possible initial cuts can become computationally intensive, since we need to run a discrete-event simulation for each evaluation of $W_{\pi(\omega)l}^\omega$. Recall though that in practice $f_\omega$ is likely to depend most strongly on  objects which are, in some sense, ``neighbours'' of $\pi(\omega)$. For instance, in Example \ref{example:intro}, the total delay for all jobs that arrive in a given time period depends strongly on the adjacent time periods. Therefore we will restrict the initial cuts to pairs which are ``nearest neighbours'' in some natural sense. Now, having laid the theoretical groundwork, the rest of this paper will focus on applications.

%% file: sec4.tex
Example \ref{example:intro} alluded to a class of scheduling problems with performance measures that have an underlying queuing structure.
In this section, we study this class of problems in greater detail.
In Table \ref{table:schedparams}, we have summarised the notation for the SAA form of the problem.

\begin{table}{
		\centering
		\rowcolors{1}{gray1}{gray2}
		\begin{tabular}{lrl}
			\textit{Sets}
			& $T$ & The set of time slices ({\footnotesize $T = [\lvert T\rvert]$}) \\
			& $S$ & The set of scenarios ({\footnotesize $S = [\lvert S\rvert]$}) \\
			& $J(s)$ & The set of jobs in scenario $s \in S$ ({\footnotesize $J(s) = [\lvert J(s)\rvert]$}) \\
			& $Y$ & The feasible region ({\footnotesize $Y \subset \mathbb{Z}_+^T$}) \\
			\textit{Parameters} 
			& $L$ & The length of a time period ({\footnotesize $L \in \mathbb{R}_+$})\\
			& $N$ & The length of the night shift ({\footnotesize $N \in \mathbb{R}_+$}) \\
			& $y_{\max}$ & The maximum number of agents available each time period \\
			& $y_{\min}$ & The minimum number of agents available each time period \\
			& $c$ & The cost per agent per time period \\
			& $d$ & The cost per unit of delay \\
			\textit{Samples}
			& $r_{sj}$ & The release time of job $j \in J(s)$ in scenario $s \in S$ \\
			& $\rho_{sj}$ & The processing time of job $j \in J(s)$ in scenario $s \in S$ \\
		\end{tabular}
		\caption{Sets and parameters for the class of job scheduling problems.}
		\label{table:schedparams}
}\end{table}

We are given a set of time periods, where by \textit{time period} $t$ we mean the half-open interval $[(t-1)L, tL) \subset \mathbb{R}$. We also have a $\left(\lvert T\rvert + 1\right)\ts{st}$ time period which we call the \textit{night shift} for practical reasons. We get to decide how many agents are allocated to each time period (\textit{during the day}). We assume there are $y_{\min}$ agents available for the night shift. We are given a set of jobs which are processed by available agents as they are released, and are scheduled on a FCFS basis. 

Note that while $j$ was used to index the objects in previous sections, we are now using $t$ to index the objects---which are time periods---and $j$ to index the jobs. Each job has a release date and a processing time which are random variables. The number of jobs may also be a random variable. Jobs consume an available agent for their duration, and
we allow overtime; in other words, an agent will always finish processing a job that starts, even if the number of agents is due to be reduced (so the number of \textit{available} agents may be negative). The goal is to minimise a positive constant multiple $d$ of the expected total delay of all the jobs; we may also incur a non-negative cost $c$ for agent time. All jobs must be finished by the end of the night shift. Since the results of this section apply equally to each scenario, for simplicity we will omit $s$ subscripts.

To get an expression for the delays we introduce the following continuous variables. For all $j \in J$, $\sigma_j \geq 0$ is the start time of job $j$.
Since the scheduling rule is fixed, the start times are completely determined by the agent levels. So there exists a mapping,
\begin{equation*}
	\sigma : \mathbb{Z}_+^T \to \prod_{j \in J} [r_j, \infty), \label{start_map}
\end{equation*}
which takes a vector $y \in \mathbb{Z}_+^T$ of agent levels to the corresponding vector, $\sigma(y) \in \mathbb{R}^J$, of start times. The delay of job $j$ is defined by 
\begin{equation}
\Delay_j(y) =
\sigma_j(y) - r_j.
\label{eq:delaydefinition}
\end{equation}
We show that the delays are non-increasing, and thus, that the problem fits the framework of (\ref{problem:main}).
We can assume that the jobs are sorted by release date.
For a fixed $y \in \mathbb{Z}_+^T$ we can compute $\sigma(y)$ recursively as follows. For $j \in J$ and $\tau \in [0,\infty)$, let
\begin{equation*}
	\kappa(\tau, j, y) = \lvert \{ j' \in J : j' < j, \sigma_{j'}(y) + \rho_{j'} > t \} \rvert 
\end{equation*}
denote the number of jobs \emph{before} $j$ which are still being processed at time $\tau$. Suppose $\sigma_1,\ldots,\sigma_{j-1}$ are known. Then $\sigma_j$ is the optimal solution to the following optimisation problem:
\begin{align}
	\min \quad & \sigma_j, \\
	\st \quad & \sigma_j \geq \max\{\sigma_{j-1}, r_j\}, \label{compute_sigma:con1} \\
	& \kappa(\sigma_j, j, y) < y_{\lceil \sigma_j/L \rceil}, \label{compute_sigma:con2}
\end{align}
where $y_t= y_{\min}$ for $t > \lvert T\rvert$. Constraint (\ref{compute_sigma:con1}) enforces the FCFS rule, and constraint (\ref{compute_sigma:con2}) ensures that job $j$ does not start unless at least one agent is available. Now we have the following lemma.
\begin{lemma}\label{lem:non_increasing_1}
	The components, $\sigma_j(y)$, of \textup{(\ref{start_map})}, are non-increasing functions. 
\end{lemma}
\begin{proof}
	It is obvious that $\sigma_1(y)$ is non-increasing. Fix $j \in J$ and suppose $\sigma_1(y),\ldots,\sigma_{j-1}(y)$ are non-increasing. Since $p_{j'}$ are constants, $\sigma_{j'}(y) + p_{j'}$ is non increasing for $j' \in \{1,\ldots,j-1\}$, which means $\kappa(\tau, j, y)$ cannot increase with $y$. Moreover, since $r_j$ is a constant, $\max\{\sigma_{j-1}(y), r_j\}$ is non-increasing. It follows that $\sigma_j(y)$ is non-increasing, which completes the proof.
\end{proof} 

Lets briefly reintroduce the $s$ subscripts. The monotonicity of $\Delay_{sj}$ is an immediate consequence of this lemma. Now note that the sum of finitely-many non-increasing functions is also non-increasing. So we can introduce Benders variables in one-to-one correspondence with any partition of the jobs; from one variable per job (\textit{per scenario}), to a single variable representing the total performance. Since the delay of a job depends strongly on the number of agents available near to its release time, it is natural to partition the jobs according to the time period they are released in. In minor abuse of notation, for each $s \in S$ and $t \in T$ let $\Delay_{st}(y)$ denote the total delay of all jobs released in time period $t \in T$ in scenario $s \in S$ given $y$, so that
\begin{equation*}
\sum_{s\in S}\sum_{t \in T}\Delay_{st}(y) = \sum_{s\in S} \displaystyle\sum_{j \in J(s)}\Delay_{sj}(y).
\end{equation*}
For each $t \in T$ we introduce a new non-negative continuous variable $\theta_{st}$ to approximate $\Delay_{st}(y)$. The precise master problem formulations are outlined in \ref{sec:8}.
At incumbent solutions to the master problem, we can compute the delays of all jobs with a discrete event simulation of the job queue for each scenario. The computational details of the DES are summarised in the appendix.
Before solving the master problem, we add initial cuts on all pairs $(t, t+1)$ of adjacent time periods.
We compute the sets $\mathcal L_{st}(y)$ which feature in a Benders cut on $(s,\,t,\, y)$ using Algorithm \ref{alg:getL} in the appendix. The idea is to expand (and contract) the set $\{t\}$ greedily until (and while) the equality is obtained (preserved).

\subsection{Integer Programming Formulation}

For benchmarking purposes, we describe a direct IP formulation. We will omit the $s$ subscripts again, but reintroducing them is easy in practice. To model the problem as an IP we need a suitable time discretisation. Let 
$\Mins = [0, \lvert T\rvert L + N)\cap \mathbb{Z}$
be the set of minutes in the time horizon, and assume that all $r_j$ and $\rho_j$ are integers. Then all jobs start and finish at discrete minutes in an optimal solution. For $j \in J$ let
\begin{equation*}
	\CanStart_j = \{r_j, \ldots, \max \Mins - \rho_j\}
\end{equation*}
denote the set of minutes that job $j$ may start, such that it will be finished by the end of the time horizon. For each minute $m \in \Mins$, define
\begin{equation*}
	\Avail_m = \{\,j \in J: r_j \leq m\}
\end{equation*} 
to be the set of jobs which have been released so far which may still start. 
For each job $j \in J$ and each minute $m \in \Mins$, we introduce a binary variable, $\sigma_{jm} \in \{0, 1\}$, with the following interpretation:
\begin{equation*}
	\sigma_{jm} = 1 \text{ if and only if job }j\text{ has started by minute } m.
\end{equation*}
As long as $Y$ is described by linear constraints, we obtain the IP formulation (\ref{ip:sched:obj} to \ref{ip:sched:con8}):

\begin{alignat}{2}
	\min \quad & \sum_{t \in T} cy_t \;+\; \sum_{j \in J}\sum_{m \in \CanStart_j}d(m - r_{j})(\sigma_{jm}-\sigma_{j,\,m-1}) \label{ip:sched:obj}
\end{alignat}
\begin{alignat}{2}
	\st \quad & y \in Y & & \label{ip:sched:con1}\\
	& \sigma_{j,\,m-1} \leq \sigma_{jm} &\qquad& \forall j \in J,\, m \in \Mins, \label{ip:sched:con2} \\
	& \sigma_{jm} \leq \sigma_{j-1,\, m} & & \forall j \in J,\,m \in \Mins, \label{ip:sched:con3} \\
	& \sigma_{jm} = 1  & & \forall j \in J,\,m \in \{\max \CanStart_j, \ldots, \max\Mins \} \label{ip:sched:con4} \\
	& \sigma_{jm} = 0 & & \forall j \in J,\, m \in \{-1,0,\ldots, r_j - 1\} \label{ip:sched:con5} \\
	& y_t \in \{y_{\min}, \ldots, y_{\max}\} & & \forall t \in T, \label{ip:sched:con6} \\
	& \sigma_{jm} \in \{0, 1\} & & \forall j \in J,\,m \in \Mins, \label{ip:sched:con7}
\end{alignat}
\begin{multline}
	\qquad\quad \sigma_{jm} - \sigma_{j, \, m - 1} \leq y_{\lceil m / L\rceil} \quad-\quad 
	\sum_{\mathclap{\substack{j' \in \Avail_{m}\setminus\{j\}, \\ m' =\, \max\{m - \rho_{j'}, r_{j'}-1\} }}}(\sigma_{j'm} - \sigma_{j'm'}) \\ \;+\; (y_{\max} - y_{\min} )(1 - \sigma_{jm} + \sigma_{j,\, m - 1}) \quad\qquad
	\forall j \in J,\, m \in \CanStart_j. \qquad \label{ip:sched:con8}
\end{multline}

The objective function measures the cost of agent hours and the total delay of all jobs.
Constraint (\ref{ip:sched:con1}) restricts the set of feasible shift schedules.
Constraint (\ref{ip:sched:con2}) ensures that job $j$ starts at exactly minute $m$ if and only if $\sigma_{jm} - \sigma_{j,\,m-1} = 1$, which allows for correct calculation of the delays. 
Constraints (\ref{ip:sched:con3}) enforces the FCFS scheduling rule.
Constraint (\ref{ip:sched:con4}) ensures that every job finishes by the end of the time horizon.
Constraint (\ref{ip:sched:con5}) ensures that every job is unprocessed at the beginning of the time horizon. Constraints (\ref{ip:sched:con6} to \ref{ip:sched:con7}) define the domains of the decision variables. Finally, constraint (\ref{ip:sched:con8}) ensures that each job does not start until there is an available agent.
The ``big M'' expression, $y_{\max} - y_{\min}$, ensures that this constraint does not prohibit overtime.
\textit{Recall too that} $y_t = y_{\min}$ \textit{for} $t > \lvert T\rvert$.
In practice, we add $s$ subscripts to the $\sigma$ variables, $\CanStart$, and $\Avail$, and modify the constraints and objective with that in mind.
We will see in Section \ref{sec:5} that this IP is intractable on realistic instances even with one scenario.

%% file: sec5.tex
In this section we will describe a \textit{Nursing Home Shift Scheduling} (NHSS) problem. This application is motivated by the rapidly increasing number of elderly people. Even though elderly people have, on average, better overall health than before, many still face disabilities or other chronic conditions that come with age. In view of this, the demand for long-term care is projected to grow substantially in the next few decades.

Health care is a prominent application area in operations research, but researchers have mostly addressed optimisation problems for hospitals, emergency services, and home care services. Studies focusing on nursing homes are scarce. The studies by \cite{LiederMoekeKooleStolletz2015}, \cite{MoekeGeerKooleBekker2016}, \cite{MoekeKoolVerkooijen2014}, \cite{vanEedenMoekeBekker2016}, \cite{BekkerEtAl2019}, and \cite{MoekeBekkerSchmidt2019} are the first and so far the only ones to address capacity planning and shift scheduling problems for nursing homes. The challenge in this context is to meet residents' requests for care or assistance within a reasonable time frame. Requests are partly scheduled in advance (such as taking medicines) and partly unscheduled (such as going to the bathroom).

Perhaps the most ambitious study to date is \cite{MoekeBekkerSchmidt2019}. Based on real-life data of both scheduled and unscheduled care in a nursing home, the authors formulate an optimisation model to determine the optimal shift scheduling for one day. To deal with the complexity and size of the problem, they approximate the optimisation problem using a Lindley recursion. The resulting MIP is then solved using standard techniques. Although the results provide important insights into the optimal staffing pattern for a nursing home during one day, the approach also has several drawbacks. Besides the fact that only an approximation of the problem is solved, the approach does not scale well: analysing larger nursing homes or extending the optimisation problem to several days seems unfeasible.

\subsection{Problem Description}

We are given a set $T$ of time periods, and a set $\Lambda$ of shift lengths in hours. Each time period has length $L$ in minutes. Shifts are allowed to start every $L$ minutes, beginning each day from time $0$. In light of practical considerations, we assume $L$ is a divisor of $60$. We let $\Gamma$ denote the set of possible shifts, which can be generated ahead of time. If $\Lambda = \{4, 8\}$ and $L = 60$, for instance, then we allow all shifts of $4$ or $8$ hours, starting every hour from time $0$. For each $t \in T$ and $\gamma \in \Gamma$ we let $\alpha_{t\gamma}$ be a binary indicator, with $\alpha_{t\gamma} = 1$ if and only if shift type $\gamma$ covers time period $t$. We will choose which shift types to use, how many to use, and schedule them. We have a minimum number $y_{\min}$ of care workers which must be available in each time period, and a maximum number $h_{\max}$ of working hours to allocate over the whole time horizon. For $y_{\max}$ we can choose a sufficiently large $y_{\max} \leq \lfloor h_{\max}/\min \Lambda\rfloor$. At the end of the day we have a night shift with $y_{\min}$ available care workers.

Scheduled and unscheduled requests arrive over the course of the day and are processed on a FCFS basis. The number and preferred start times of scheduled requests are known in advance, but the the unscheduled requests are stochastic, as is the duration of all requests.
Requests that we fail to service by the end of the day are serviced by the night staff instead. We incur no cost for staffing hours (although we could), and the goal is to minimise the expected total delay of all requests.

To see that the NHSS problem is a special case of the general framework, note that as usual, the start times of the requests are determined by the number of care workers available in each time period. These levels are, in turn, uniquely determined by the shift schedule. The feasible region $Y$ will be the set of vectors of care worker levels that constitute a legal shift schedule, which is easily described using linear constraints. The staffing hours are the resources, the time periods are the objects, and the waiting time of patients corresponds to the delay of the jobs. The MP of the LBBD formulation is detailed in \ref{sec:8}.

%

\subsection{Instance Generation}

For our computational experiments, we utilise fictional data based on the statistical analysis of nursing home operations provided by \cite{vanEedenMoekeBekker2016} and \cite{MoekeBekkerSchmidt2019}. A work day consists of $16$ hours from $7$ AM to $11$ PM, with the night shift spanning $11$ PM to $7$ AM. In view of safety regulations, as well as standard practice, we require a minimum of $y_{\min}=2$ care workers to be scheduled at all times, with a (very conservative) maximum of $y_{\max}=20$ care workers at any given time. We have a maximum of $h_{\max} = 80$ working hours. There are two shift types: short shifts have a duration of 4 hours, while long shifts have a duration of 8 hours. A shift can start at every full hour. Thus we take $\Lambda = \{4, 8\}$ as the shift types, and divide the day into intervals of length $L = 60$ minutes. The dataset contains $224$ scheduled requests which are available as a csv file in the associated Github. The unscheduled requests can either be of the short type or the long type, and each request has an $90$\% chance of being short. The duration of each short requests is sampled from an exponential distribution with mean $1.89$ minutes, while long requests use a mean of $9.28$ minutes. Unscheduled requests follow a Poisson process with rate $\lambda=20$. We will refer to instances generated according to this scheme as Type A instances.

We will see shortly that Type A instances are too hard for the IP formulation. So to more clearly compare the approaches, we will also consider a second family of parameters. Type B instances are characterised by a $4$ hour work day and a $2$ hour night shift. We require a minimum of $y_{\min} = 1$ care workers and allow at most $y_{\max}=10$ workers at a time. We have a maximum of $h_{\max} = 15$ working hours and allow shifts of $1$ or $2$ hours. The scheduled and unscheduled requests are generated identically to Type A instances, but we divide all release times by four and discard all but every fourth request. To accommodate the IP formulation, all data is rounded to positive integer values, but this is not required by the LBBD approach.

\subsection{Benchmarking}

All implementations in this paper were programmed in Python 3.8.8 via the Anaconda distribution (4.9.2). Master problems were solved using Gurobi 9.1.2 (\cite{Gurobi}), and discrete event simulations were coded from scratch. Jobs were run in parallel on a computing cluster operating $2.5$GHz CPUs. Each job was allocated at most $5$GB of memory and a single thread. We were typically able to achieve faster solution times using a typical desktop PC but opted for the cluster to ensure our benchmarking was rigorous.

We consider Type $B$ instances first. Our aim is to compare the the LBBD and IP formulations with respect to runtime. To that end we generated $20$ Type B instances with a single scenario each. We attempted to solve each instance with three approaches:
\begin{equation*}
\begin{tabular}{rl}
	\textsf{BD} & Logic-Based Benders Decomposition, \\
	\textsf{IP} & Direct integer program, and \\
	\textsf{IP}\textit{f} & Direct integer program with the optimal shift schedule fixed.
\end{tabular}
\centering
\end{equation*}

In the case of \textsf{IP}\textit{f} the $y$ variables of the model are fixed to the optimal solution found using the other methods. In other words, \textsf{IP}\textit{f} solves the pure job scheduling problem associated with an optimal shift schedule, while \textsf{IP} denotes the full integer programming formulation. The times required to prove optimality are documented in Table~\ref{tbl:benchmarking}, where we see that LBBD  outperforms the IP model by several orders of magnitude on every instance. Indeed, even solving the pure job scheduling problem associated with the optimal shift schedule is significantly slower than the full LBBD formulation.

\begin{table}[h!]
	\rowcolors{1}{gray1}{gray2}
\begin{tabular}{rcccc}
	ID & Objective ($\prime$) & \textsf{BD} & \textsf{IP} & \textsf{IP}\textit{f} \\
	0 & 64 & 0.0092 & 589.9901 & 30.9762 \\
	1 & 22 & 0.0124 & 137.8063 & 41.0888 \\
	2 & 18 & 0.0044 & 60.6626 & 37.8141 \\
	3 & 117 & 1.1238 & 1930.5293 & 29.4467 \\
	4 & 147 & 0.017 & 2451.8047 & 31.4453 \\
	5 & 13 & 0.0097 & 87.3941 & 50.2346 \\
	6 & 18 & 0.0126 & 93.3177 & 46.6182 \\
	7 & 42 & 0.0115 & 444.6288 & 39.4959 \\
	8 & 19 & 0.0148 & 81.2079 & 36.8963 \\
	9 & 24 & 0.0194 & 519.7442 & 56.1756
\end{tabular}\quad
\begin{tabular}{rcccc}
	ID & Obj. & \textsf{BD} & \textsf{IP} & \textsf{IP}\textit{f} \\
	10 & 21 & 0.0091 & 97.4927 & 28.6019 \\
	11 & 26 & 0.0045 & 322.773 & 46.9173 \\
	12 & 35 & 0.013 & 418.8719 & 53.5033 \\
	13 & 9 & 0.0115 & 56.6321 & 40.2135 \\
	14 & 38 & 0.009 & 510.0433 & 29.0088 \\
	15 & 124 & 0.0112 & 1209.5284 & 30.9684 \\
	16 & 41 & 0.0099 & 567.0863 & 31.241 \\
	17 & 74 & 0.0115 & 596.8426 & 34.6893 \\
	18 & 18 & 0.0084 & 72.6107 & 40.5433 \\
	19 & 131 & 0.0146 & 2568.7499 & 32.0914
\end{tabular}
\centering
\caption{Solution times for Benders Decomposition, the IP model, and the IP model with an optimal shift schedule fixed, on $20$ Type B instances of the NHSS problem.}
\label{tbl:benchmarking}
\end{table}

\subsection{Further Results}

Based on the results for the Type B instances, it is no surprise that the IP model is intractable for Type A instances. In this section we consider not only Type A instances, but multiple scenarios. We generated 10 instances each with $\lvert S\rvert = 1, 10$, and $100$ scenarios. The solution times of the $\textsf{BD}$ method are found in Table~\ref{tbl:nursingscenariostimes} together with more detailed solver information for the $\lvert S\rvert = 100$ instances.

\begin{table}[h!]
	\centering
	\rowcolors{1}{gray1}{gray2}
	\begin{tabular}{rccccccc}
		& \multicolumn{3}{c}{Total Runtimes (s)} & \multicolumn{4}{c}{Further $\lvert S\rvert = 100$ results} \\
		ID & $\lvert S\rvert = 1$ & $\lvert S\rvert = 10$ & $\lvert S\rvert = 100$ & Objective ($\prime$) & Initials Time (s) & Callback Time (s) & Cuts \\
		0 & 3.63 & 7.72 & 57.77 & 750.98 & 17.28 & 6.60 & 854 \\
		1 & 0.25 & 2.13 & 36.17 & 736.71 & 16.82 & 0.36 & 412 \\
		2 & 2.57 & 5.31 & 40.81 & 758.05 & 17.28 & 0.27 & 336 \\
		3 & 1.54 & 4.85 & 69.15 & 764.95 & 17.28 & 7.27 & 892 \\
		4 & 0.24 & 4.94 & 56.83 & 775.56 & 17.52 & 5.88 & 982 \\
		5 & 1.58 & 6.94 & 44.81 & 782.73 & 17.51 & 3.54 & 654 \\
		6 & 2.55 & 7.95 & 40.16 & 785.69 & 17.45 & 0.18 & 308 \\
		7 & 2.30 & 5.67 & 62.92 & 764.52 & 17.49 & 4.73 & 880 \\
		8 & 1.92 & 4.78 & 56.48 & 791.05 & 17.55 & 2.76 & 844 \\
		9 & 1.47 & 4.28 & 78.14 & 815.97 & 17.69 & 8.85  & 1774
	\end{tabular}
\caption{Solution times for LBBD on Type A instances of the NHSS problem with $\lvert S\rvert = 1, 10$, and $100$ scenarios, with additional information about the $\lvert S\rvert = 100$ case.}
\label{tbl:nursingscenariostimes}
\end{table}

We see that even with $100$ scenarios, all instances can be solved in a reasonable amount of time. Typically around one third of the solution time was spent generating initial cuts. We will study the importance of the initial cuts in more detail in the context of the ACCA problem, where we will solve identical instances with \textit{and} without them. Note also that despite the large number of feasible solutions, fewer than $1000$ Benders cuts were needed to prove optimality for all but one instance. We will also look at this phenomenon more closely in the context of the ACCA problem.

For now, we are interested in the accuracy of the sample average approximation. To that end, we cross validated all ten of the $\lvert S\rvert = 100$ instances. In other words, we calculated the optimal objective value of each optimal solution by passing it to the simulation for every other set of scenarios. The percentage differences between the optimal objective value for each set of ``training" scenarios (that is, those scenarios we optimised over to obtain a given solution), and value of that solution using the various ``test" scenarios, are recorded in Table \ref{tbl:nhss:crossvalidate}. The very small values off the diagonal show that we have successfully captured the inherent variance of the problem.

\begin{table}
\rowcolors{1}{gray1}{gray2}
\begin{tabular}{rcccccccccc}
	& \multicolumn{10}{c}{Test Instance} \\
	Training Instance & 0 & 1 & 2 & 3 & 4 & 5 & 6 & 7 & 8 & 9 \\
	0 & 0 & 0.038 & 0.033 & 0 & 0.018 & 0.01 & 0.024 & 0.01 & 0.003 & 0.006 \\
	1 & 0.021 & 0 & 0.018 & 0.011 & 0.009 & 0.016 & 0.026 & 0.029 & 0.011 & 0.006 \\
	2 & 0.006 & 0.044 & 0 & 0.002 & 0 & 0 & 0.003 & 0.001 & 0.03 & 0 \\
	3 & 0 & 0.038 & 0.033 & 0 & 0.018 & 0.01 & 0.024 & 0.01 & 0.003 & 0.006 \\
	4 & 0.006 & 0.044 & 0 & 0.002 & 0 & 0 & 0.003 & 0.001 & 0.03 & 0 \\
	5 & 0.006 & 0.044 & 0 & 0.002 & 0 & 0 & 0.003 & 0.001 & 0.03 & 0 \\
	6 & 0.005 & 0.04 & 0.003 & 0.005 & 0.004 & 0.004 & 0 & 0 & 0.022 & 0.001 \\
	7 & 0.005 & 0.04 & 0.003 & 0.005 & 0.004 & 0.004 & 0 & 0 & 0.022 & 0.001 \\
	8 & 0.009 & 0.045 & 0.046 & 0.014 & 0.025 & 0.029 & 0.031 & 0.018 & 0 & 0.012 \\
	9 & 0.006 & 0.044 & 0 & 0.002 & 0 & 0 & 0.003 & 0.001 & 0.03 & 0 \\
\end{tabular}
\centering
\caption{Cross validation of the optimal objective value on $10$ distinct Type A instances Type A instances of the NHSS problem with $\lvert S\rvert = 100$. The entry in row $i$ and column $j$ is $\lvert a - b\rvert/a$, where $a$ is the objective value for instance $i$, and $b$ is the value of the optimal solution to instance $i$ using the scenarios from instance $j$ (rounded to three places). }
\label{tbl:nhss:crossvalidate}
\end{table}

\begin{figure}[h!]
	\begin{center}
	\includegraphics[width=.3\linewidth]{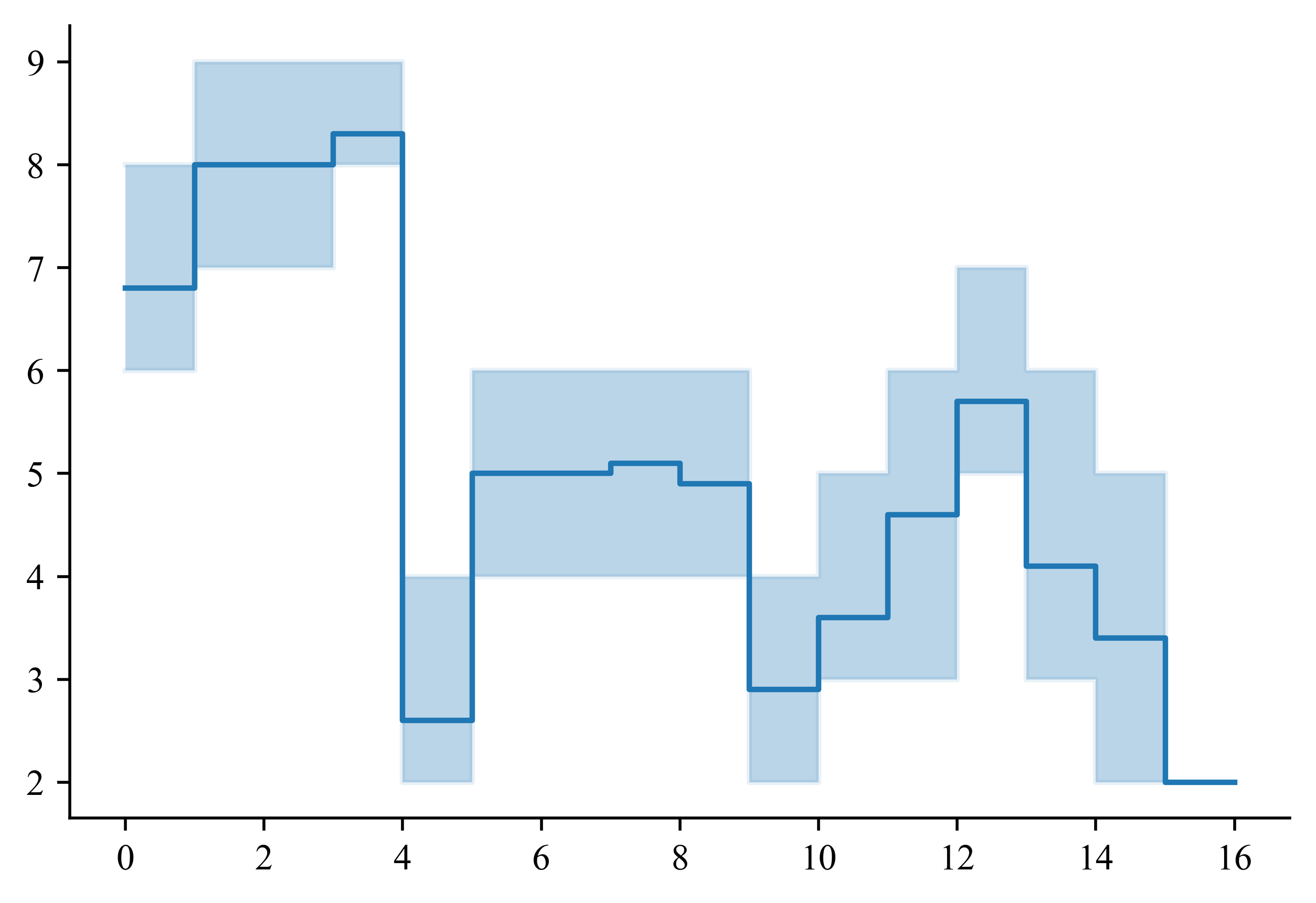}
	\includegraphics[width=.3\linewidth]{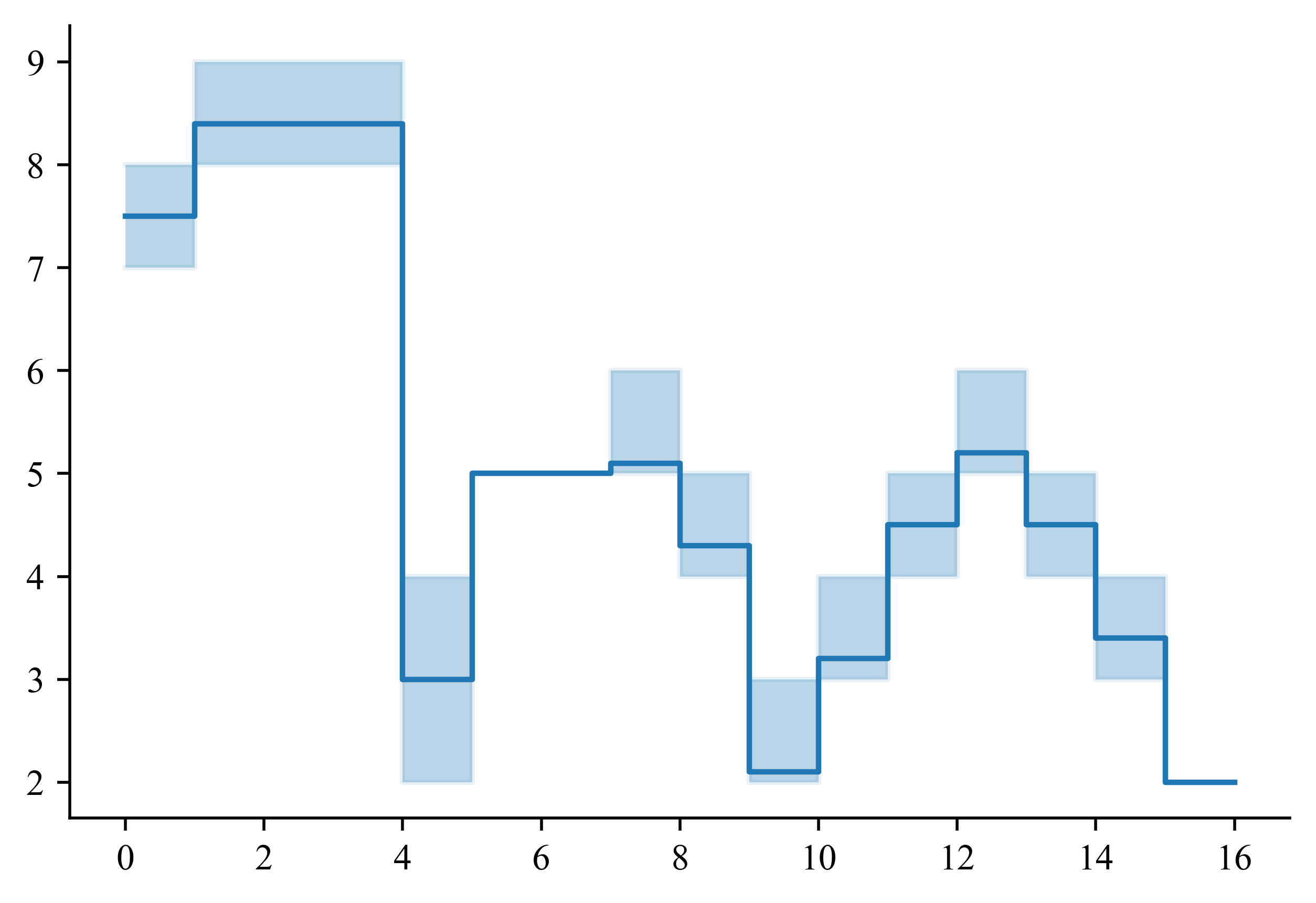}
	\includegraphics[width=.3\linewidth]{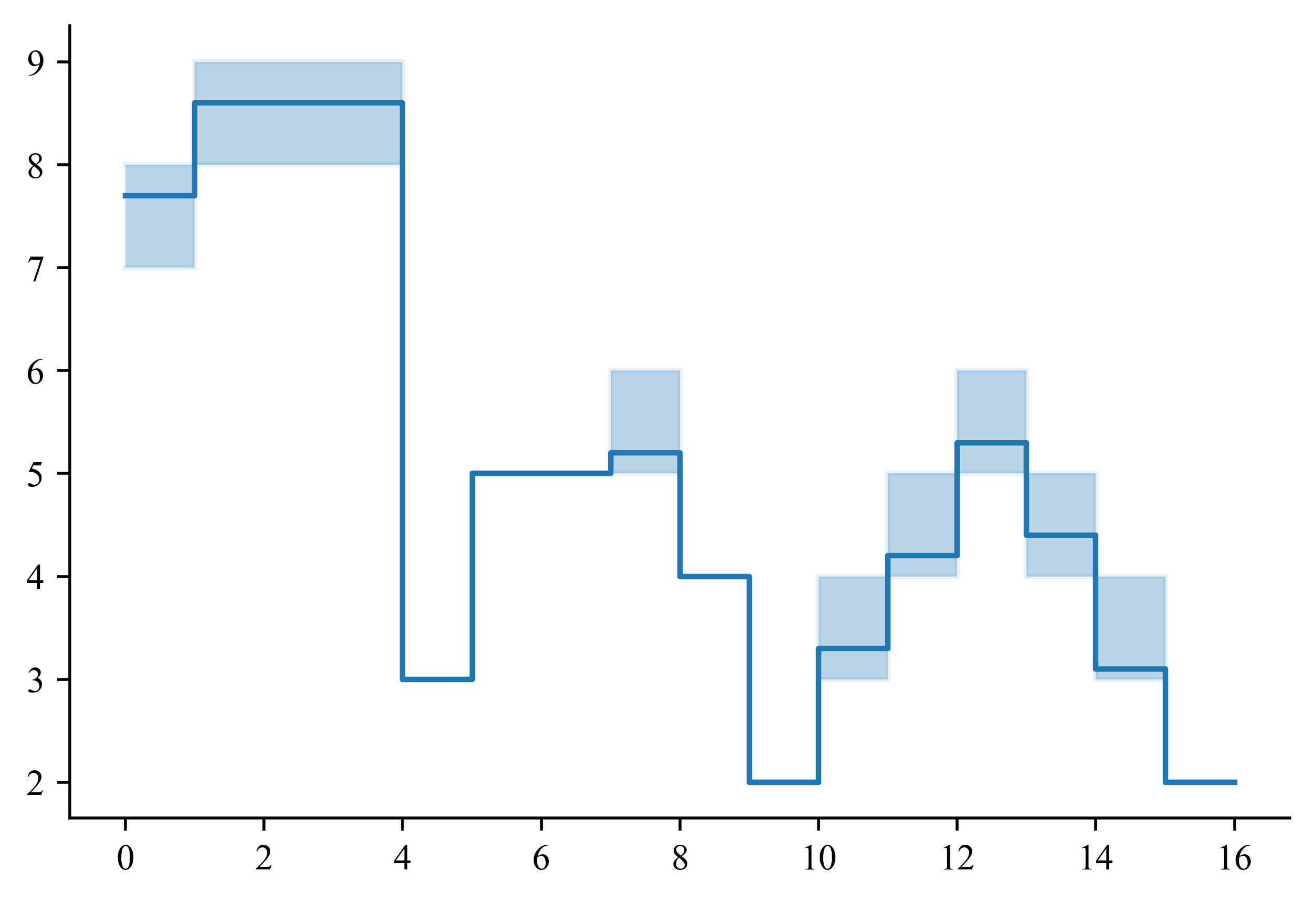}
	\end{center}
	\caption{Variance in the optimal shift schedule for $\lvert S\rvert = 1$ (left), $=10$ (middle) and $=100$ (right) scenarios.}
	\label{fig:nhss:stability}
\end{figure}

Lastly we examine the stability of the optimal solution. The range and averages of the optimal number of care workers for each time period are illustrated in Figure \ref{fig:nhss:stability} for $\lvert S\rvert = 1, 10$, and $100$ instances.
Clearly one scenario is not enough to produce a stable solution. And while $\lvert S\rvert = 10$ scenarios is a significant improvement, with $100$ scenarios there is little variance in the optimal shift schedule. In other words, it does not matter \textit{which} $100$ scenarios we optimise over, as long as we optimise over enough of them.

%% file: sec6.tex
In this section we study another concrete example of the scheduling problem introduced in Section \ref{sec:4}, an \textit{Airport Check-in Counter Allocation} (ACCA) problem. This application is motivated by the increase in passenger air traffic witnessed over the past few decades. Modern airports tend to service a large number of distinct airlines, flights, and passengers. Two major choke points on passenger flow within an airport are security checkpoints, and check-in terminals. Although self-service check-in terminals have grown in popularity, traditional staffed checkpoints are likely to remain necessary for several reasons, including security concerns, baggage handling, and personal preference. In any case, the model we develop is applicable to a wide variety of multi server queuing problems.

An early integer programming model was proposed by \cite{DijkSluis2006} to minimise the number of check-in counters required to service all flights. Simulation modelling is then used to evaluate the performance of the solutions produced by the approximate optimisation model. We will test our model on a numerical example derived from this paper.

\cite{BrunoGenovese2010} propose integer programming models for determining the number of check-in counters to open. \cite{AraujoRepolho2015} extend their work by introducing a service level constraint, which caps the proportion of passengers still in the queue at the end of each time period. \cite{BrunoDiglioGenovesePicollo2019} extend the model further by incorporating a staff scheduling component. \cite{LalitaMannaMurthy2020} propose an integer programming model that attempts to enforce FCFS queue discipline by ensuring that sufficient check-in capacity is available. These models account for the stochasticity of the problem by including estimates of the passenger arrival rates and service times as model parameters. Thus they fail to account for any variance in the parameters. Furthermore, they embed only an approximation of the queuing structure into the model. 

\subsection{Problem Description}

We are given a set $T$ of time periods, a set $F$ of flights, and a set $J$ of passengers.
Each time period has length $L$, and each passenger is scheduled to board one of the flights. 
The passengers arrive over the course of the time horizon, and must pass through the check-in terminal. 
Passengers are served on a FCFS basis at open check-in counters, and we get to decide how many check-in counters to open in each time period. 
We can open up to $y_{\max}$ check-in counters at a time, and allow as few as $y_{\min}=0$ to be open at any time.
We incur a fixed cost of $c=D$ for opening a check-in counter for one time period, as well as a cost of $Q$ per time period for the total delay (that is, queuing time). As usual we allow overtime; in other words, we will finish serving all passengers even if the number of open check-in counters is due to reduce.

For each flight $f \in F$ we have a deadline $d_f$, which is the latest time we are allowed to finish servicing any passenger in order for them to be able to make their flight. We have no night shift ($N = 0$) and instead incur a penalty of $\lvert T\rvert\times L$ for each passenger that we fail to serve by their departure time. This modification clearly preserves the monotonicity of the performance measures. To ensure that the passenger deadlines can be met in every scenario, we will introduce two families of constraints.
For each $f \in F$ let $J_f \subset J$ denote the set of passengers for that flight.
For each $t \in T$ define
\begin{equation}
	B_t = \max_{s \in S}\sum_{\substack{f \in F \\ :d_f < Lt}} \sum_{\substack{j \in J_f}}\rho_{sj} \qquad
	A_t = \max_{s \in S} \sum_{\substack{j \in J \\ :Lt \leq r_{sj}}}\rho_{sj}.
\end{equation}
$B_t$ is the maximum total service time of passengers for flights which depart during or before time period $t$, while $A_t$ is the maximum total service time of all passengers who arrive during or after time period $t$. Service level constraints of the form 
\begin{equation*}
	\sum_{t' = t}^{\lvert T\rvert}L y_{t'} \geq A_t, \quad \text{and} \quad \sum_{t' = 1}^t Ly_{t'} \geq B_t \quad \text{for }t \in T
\end{equation*}
ensure that sufficient check-in capacity is made available to clear all passengers in every scenario, as long as $y_{\max}$ is sufficiently high. The MP of the LBBD formulation is detailed in \ref{sec:8}.

\subsection{Instance Generation}

To generate instances for the ACCA problem, we used the numerical example from \cite{DijkSluis2006}. We have $21$ time periods, each of which are $L=30$ minutes long. There are $10$ flights, and $2160$ passengers. The number of passengers for each flight is given in (\ref{table:flights}).

\begingroup\small
\begin{equation}
	\rowcolors{1}{gray1}{gray2}
	\begin{tabular}{lcccccccccc}
		Flight & 1 & 2 & 3 & 4 & 5 & 6 & 7 & 8 & 9 & 10 \\
		Passengers & 150 & 210 & 240 & 180 & 270 & 150 & 210 & 300 & 180 & 270 \\
		Starting time & 0 & 2 & 4 & 4 & 6 & 8 & 10 & 12 & 12 & 14 \\
	\end{tabular}
	\label{table:flights}
\end{equation}
\endgroup

\smallskip	
For each flight we have an ``arrival period'' which consists of seven sequential time periods. The first time of each arrival period is given in the second row of (\ref{table:flights}). All passengers for a flight arrive in the corresponding arrival period. The distribution of passenger arrivals over an arrival period is given in (\ref{table:arrivals}).

\begingroup\small
\begin{equation}
	\rowcolors{1}{gray1}{gray2}
	\begin{tabular}{lccccccc}
		Time period & 1 & 2 & 3 & 4 & 5 & 6 & 7 \\
		\% of passengers & 5 & 10 & 20 & 30 & 20 & 15 & 0 
	\end{tabular}
	\label{table:arrivals}
\end{equation}
\endgroup

The deadline for each flight is the end of its arrival period (note that no passengers arrive in the final time period of their arrival period). To generate a set of scenarios, for each $s \in S$ and each $j \in P$, the arrival time of passenger $j$ in scenario $s$ is given by
\begin{equation*}
	a_{sj} \sim L \times (\mathfrak W + \textrm{Start}_{f_j}) + \mathfrak U(\{1,\ldots,L\}),
\end{equation*}
where $\mathfrak U$ denotes the discrete uniform distribution, $f_j$ is the flight for passenger $j$, $\textrm{Start}_{f_j}$ is the starting time period for flight $f_j$, and $\mathfrak W$ is the discrete random variable with the distribution given by (\ref{table:arrivals}). Service times were sampled from an exponential distribution with mean $2$. In the objective function, we incur a fixed cost of $D = 40$ for opening a check-in counter for one time period. In each time period we can open up to $y_{\max} = 20$ check-in counters.

\subsection{Computational Results}
Since all instances of this problem are larger and more difficult to solve than the NHHS problem studied in Section \ref{sec:5}, we do not bother with the IP formulation. Table \ref{tbl:airportlbbd} contains solution information for the LBBD approach on 
$10$ sets of $\lvert S\rvert = 100$ scenarios. We used a queuing cost of $Q=40$ as is standard in the literature. We see that using LBBD we are able to solve all instances to optimality in a few minutes. The time spent running simulations to generate Benders cuts in the callback routine sits at $10$ to $20\%$ of the total solver time (that is, total time minus the time spent generating initial cuts). In the vicinity of $470000$ simulations were run per solve, which includes initial cut generation. What immediately stands out is how much of the total time was spent generating initial cuts. 

\begin{table}[h!]
	\rowcolors{1}{gray1}{gray2}
	\begin{tabular}{rcccccc}
		ID & Objective & Total Time (s) & Initials Time  & Callback Time & Cuts Added & Simulations Run \\
		0 & 8442.478 & 249.073 & 219.066 & 3.895 & 1243 & 470532  \\
		1 & 8420.75 & 276.23 & 225.721 &  10.983 & 3196 & 471148 \\
		2 & 8427.365 & 264.334 & 220.455 &  13.705 & 4541 & 471195 \\
		3 & 8418.309 & 245.992 & 221.397 &  5.515 & 1791 & 470236 \\
		4 & 8432.798 & 261.041 & 229.252 & 9.024 & 2735 & 468138\\
		5 & 8463.651 & 320.047 & 227.158 & 19.066 & 4268 & 469189 \\
		6 & 8444.105 & 328.748 & 235.592 & 11.096 & 4140 & 470767 \\
		7 & 8459.151 & 316.923 & 244.668 & 8.714 & 2838 & 468924 \\
		8 & 8465.96 & 303.009 & 233.585 & 9.812 & 3457 & 471168 \\
		9 & 8438.479 & 243.797 & 211.588  & 8.674 & 3259 & 471212 \\
	\end{tabular}
	\centering
	\caption{Solution information for the ACCA problem with a queuing cost of $Q = 40$.}
	\label{tbl:airportlbbd}
\end{table}

\begin{figure}
	\begin{center}
		\includegraphics[width = \linewidth]{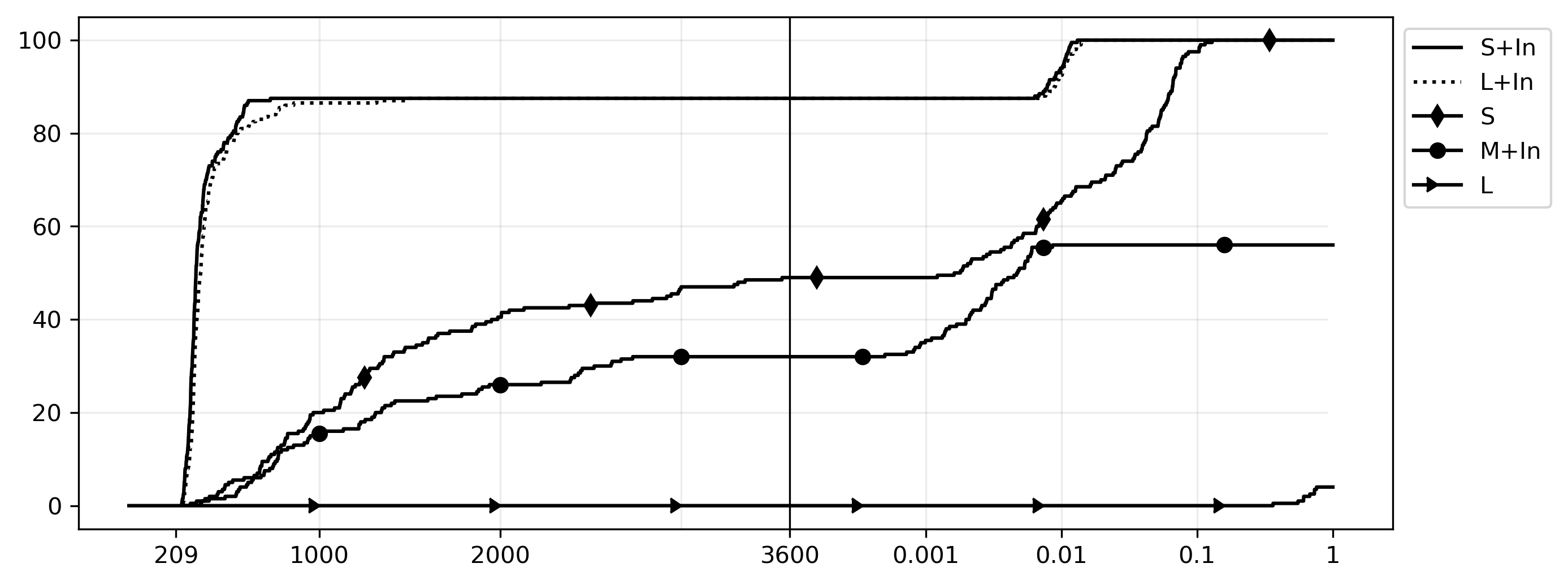}
	\end{center}
		\caption{Performance profiles for five approaches. The horizontal axis represents runtime, and the optimality gap obtained. To the left of $3600$ seconds, the vertical axis indicates the proportion of instances solved to optimality within a given time. To the right of $3600$ seconds (which also represents an optimality gap of $0$) the vertical axis indicates the proportion of instances which obtained a gap at least this small by the time limit.}
	\label{figure:airport:performance}
\end{figure}

In order to compare the relative strength of the various Benders cuts and initial cuts, we generated another $25$ sets of $\lvert S\rvert = 100$ scenarios. We considered queuing costs of $Q \in \{5, 10, 15, \ldots, 40\}$ for a total of $200$ unique instances. We tried to solve all $200$ instances with each of the four Benders cuts: the no-good cut, (\ref{cut:1}) or \textsf{NG}, the monotonic cut, (\ref{cut:2}) or \textsf{M}, the local cut, (\ref{cut:3}) or \textsf{L}, and the strengthened cut, (\ref{cut:4}) or \textsf{S}. We did this both with and without initial cuts, \textsf{+In}, for a total of eight approaches. Each run was allowed $3600$ seconds of computational time under the same conditions as in Section \ref{sec:5}. $\textsf{NG}, \textsf{NG+In}$, and $\textsf{M}$ were all unable to reach the time limit without exceeding the memory limit for every instance. Performance profiles for the five remaining methods are provided in Figure \ref{figure:airport:performance}. The $209$ second mark represents the minimum amount of time taken to generate the $84000$ initial cuts.

%

The performance profiles show that both $\textsf{L+In}$ and $\textsf{S+In}$ were able to solve the majority of instances to optimality within $1000$ seconds, and achieve a small optimality gap for the rest. Although $\textsf{L+In}$ and $\textsf{S+In}$ are comparable, the relative strength of $\textsf{S}$ is made clear by the disparity between $\textsf{S}$ and $\textsf{L}$. Without initial cuts, even $\textsf{L}$ exceeded the memory limit on all but four instances, while $\textsf{S}$ still solved almost half of the instances to optimality within the allowed time. This demonstrates the strength both of $\textsf{S}$ and the initial cuts.

Lastly, we are also interested in the number of simulations run, as a proportion of the feasible solutions. The service level constraints make it challenging to compute the precise number of feasible solutions to a given instance. For the sake of argument, consider the first (ID $0$) set of scenarios from the performance profile above. First we optimised the master problem with the new objective function, $\sum_{t \in T}y_t$. In doing so we found that to satisfy the service level constraints, we must open at least $153$ check in counters over the $21$ time periods. We then added a constraint of the form $\sum_{t \in T}y_t = 153$, and re-optimized the model with the objective, $\min \max_{t \in T}y_t$. In doing so we found that
\begin{align*}
	\bar y = (1, 0, 0, 0, 9, 9, 9, 9, 9, 9, 8, 9, 9, 9, 9, 9, 9, 9, 9, 9, 9)
\end{align*}
is a feasible solution to this instance. Since $\bar y$ is feasible, every solution which is component-wise larger than $\bar y$ is also feasible. Therefore 
\begin{equation*}
	\prod_{t \in T}\left(21 - \bar y_t\right) \approx 3.727\times 10^{23}
\end{equation*} 
is a lower bound on the total number of feasible solutions for this instance. It is a \emph{lower} bound because we have not accounted for feasible solutions which we can obtain by increasing one level and decreasing another. In any case, among the $200$ instances solved with $\textsf{S+In}$ detailed above, the \emph{maximum} number of simulations run was $603948$, which is a tiny proportion of the feasible solutions. This demonstrates the soundness of simulating incumbent solutions as a methodology.

%% file: sec7.tex
In this paper we have proposed a novel approach to solving certain complicated optimisation problems. The idea is to simulate the performance of each incumbent solution to an appropriately relaxed model. The simulation data then informs the trajectory of the optimisation model itself via logic based Benders cuts. Strong Benders cuts ensure that we converge on an exact solution to the original problem, and that the number of simulations we actually need to run is relatively small. We tested the approach on a class of stochastic resource allocation problems with monotonic performance measures. In particular we considered two concrete special cases: a nursing home shift scheduling problem, and an airport check-in counter allocation problem. 
Previous approaches failed to account for the the non-trivial variance of the underlying distributions. Optimising over a set of samples large enough to capture this variance greatly increases the difficulty of the problem. And even for few scenarios, trying to embed the queuing structure into an IP leads to an intractable model. We have seen, however, that using LBBD we are able to solve these problems to optimality with as many as $100$ scenarios in a reasonable amount of time. We are able to incorporate the queuing structure exactly at the level of scenarios, and obtain stable, accurate solutions to the stochastic problem.

We should not overlook the fact that our results were achieved despite the imposition of significant limitations on the computational resources. First of all, all computations were done on a single thread. But since the simulation is independent of the optimisation model, there is significant scope for parallelisation; especially of the expensive initial cuts. Moreover the simulations themselves were coded from scratch in Python. Further speed ups could be achieved using a faster compiled language. This improvement is not accessible to the pure IP model, since the solvers typically already use extremely efficient compiled code.

Unlike its classical counterpart, logic based Benders decomposition does not yet enjoy the status of a ``standard technique.'' We believe this paper constitutes important progress in that direction. We have seen that it is possible to feed simulation data into an optimisation problem using Benders cuts. We believe there is significant scope to apply this methodology to other optimisation problems, where incorporating certain complicated structures would render a direct approach nonviable. In particular we hope practitioners will apply our theoretical development to other resource allocation problems with monotonic performance measures. We conclude this paper by discussing some potential avenues for future research.

The LBBD approach does not depend on the structure of the feasible region of MP. So as not to distract the reader from our primary message—--the estimation of complicated performance
measures via Benders cuts—--we have omitted some potentially interesting details. Consider
the NHSS problem. It may be desirable to limit the number of four hour shifts used. We
could also consider a finer subdivision of the time horizon, say into 15 or 30 minute intervals.
In fact, we were able to solve instances with the time horizon extended to fourteen days, but this is
only interesting if the scheduled request patterns are substantially different on different days.
In the context of the ACCA problem, it may be desirable to impose additional constraints
pertaining to the structure of the queue, such as maximum queue lengths or queuing times.
During a simulation it is easy to detect if qualitative constraints like these are violated by
the current solution. We can eliminate these solutions using feasibility cuts. We can even do
much better than a no-good cut, since we know that at least one more check-in counter will
need to be opened in the set $\mathcal L$ prior to the problematic time period. It would be interesting
to pursue this in more detail.

For ease of exposition, we have considered only one class of resources and one family of performance measures at a time. Extending the approach to multiple classes is straightforward.

The sample-average approximation approach gives equal importance to each scenario. It may also be interesting to optimise other stochastic measures, such as the Conditional Value-at-Risk.


The Benders cuts and initial cuts derived in this paper are valid because the performance measures are non-increasing. The question arises: what if the cuts are only approximately non-increasing? In other words, what if adding more resources can \textit{sometimes} make the solution worse. In that case, some of the Benders cuts will be too tight and we may cut off an optimal solution. But does still constitute a good heuristic?

%% file: Acknowledgements.tex
Mitchell Harris is supported by an Australian Government Research Training Program (RTP) scholarship.
Marijn Jansen and Thomas Taimre were partly funded by Australian Research Council (ARC) Discovery Project DP180101602.
We thank Marko Boon for his help in starting up this research project.
We are grateful to Ren\'e Bekker and Dennis Moeke for sharing a fictional data set.

%% file: sec8.tex

We generate the sets $\mathcal L_{st}(y')$ that feature in Benders cuts with Algorithm \ref{alg:getL}.

\begin{algorithm}
	\caption{Calculating $\mathcal L_{st}(y')$}\label{alg:getL}
	\begin{algorithmic}[1]
		\medskip
		\State \textbf{Input:} $s,\, t,\, y'$
		\State Put $\mathcal L = \{t\}$
		\While{$\Delay_{st}(\Delta(\mathcal L, y')) < \Delay_{st}(y')$}
		\State $\mathcal L \longleftarrow \{\max(\min \mathcal L - 1, 1), \ldots, \min(\max \mathcal L + 1, \lvert T\rvert)\}$
		\EndWhile
		\While{$\Delay_{st}(\Delta(\mathcal L \setminus \min\mathcal L, y')) = \Delay_{st}(y')$}
		\State $\mathcal L \longleftarrow \mathcal L \setminus \min\mathcal L$
		\EndWhile
		\While{$\Delay_{st}(\Delta(\mathcal L \setminus \max\mathcal L, y')) = \Delay_{st}(y')$}
		\State $\mathcal L \longleftarrow \mathcal L \setminus \max\mathcal L$
		\EndWhile
		\State \textbf{Return:} $\mathcal L_{st}(y') = \mathcal L$
	\end{algorithmic}
\end{algorithm}

The master problem of the NHSS problem is (\ref{mp:nhss:start} to \ref{mp:nhss:end}):
\begin{alignat}{2}
\min\quad & \sum_{s \in S}\sum_{t \in T}\theta_{st}/S & \quad & \label{mp:nhss:start} \\
\st\quad & \sum_{\xi = y_{\min}}^{y_{\max}}z_{t\xi} = 1 & & \forall t \in T, \label{mp:nhss:yz1} \\
& \sum_{\xi = y_{\min}}^{y_{\max}}\xi z_{t\xi} = y_t & & \forall t \in T, \label{mp:nhss:yz2} \\
& \sum_{t \in T}L\cdot y_t \leq 60\cdot h_{\max} \\
& y_t = \sum_{\gamma \in \Gamma}\alpha_{t\gamma}x_\gamma \\
& x_\gamma \in \mathbb{Z}_+ & & \forall \gamma \in \Gamma,\\
& y_t \in \{y_{\min}, \ldots, y_{\max}\} & & \forall t \in T, \label{mp:nhss:y} \\
& z_{t\xi} \in \{0, 1\} & & \forall \xi \in \{y_{\min}, \ldots, y_{\max}\},\, t \in T, \label{mp:nhss:z} \\
& \theta_{st} \geq 0 & & \forall s \in S,\, t \in T. \label{mp:nhss:end}
\end{alignat}

\bigskip
The master problem of the ACCA problem is (\ref{mp:acca:start} to \ref{mp:acca:end}):
\begin{alignat}{2}
	\min\quad & \sum_{t \in T}D\cdot y_t + \sum_{s \in S}\sum_{t \in T}(Q / L)\theta_{st}/S & \quad & \label{mp:acca:start} \\
	\st\quad & \textup{(\ref{mp:nhss:yz1}), (\ref{mp:nhss:yz2}), (\ref{mp:nhss:y}), (\ref{mp:nhss:z}), (\ref{mp:nhss:end})} \\
	& \sum_{t' \in \{t,\ldots, \lvert T\rvert\}}L\cdot y_{t'} \geq A_t & & \forall t \in T, \\
	& \sum_{t' \in \{1,\ldots,t\}}L\cdot y_{t'} \geq B_t & & \forall t \in T, \label{mp:acca:end}
\end{alignat}

\bigskip
The discrete-event simulation is summarised in Algorithm \ref{alg:scheduling:main}.
\begin{algorithm}[h!]
	\begin{algorithmic}[1]
		\caption{Discrete Event Simulation}\label{alg:scheduling:main}
		\State \textbf{Input:} $s$, $\mathcal L \subseteq \{1,\ldots, n\}$, $\{(t, y_t'): t \in \mathcal L\}$
		\State $y_t'\longleftarrow y_{\max}$ for $t \in T\setminus\mathcal L$ 
		\State $f_t \longleftarrow 0$ for $t \in T$ // \textit{Performance measures}
		\State $\mathcal Q \longleftarrow (0, \ldots, 0)$ where $\len(\mathcal Q) = y_1'$
		\State $t' \longleftarrow 0$ // \textit{Current time period }
		\For{$j \in J$ }
		\State // \textit{If there are no staff available, or the next job has not arrived yet, advance}
		\While{$\len(Q) = 0$ \textbf{or} $\max\left\{ r_{sj},\, \mathcal Q_0 \right\} \geq (t' + 1)L$}
		\If{$y_{t'+1}' > y_{t'}'$} // \textit{If we increase the resources}
		\State Add $(t' + 1)L$ to $\mathcal Q$
		\EndIf
		\If{$y_{t'+1}' < y_{t'}'$} // \textit{If we reduce the resources}
		\State Remove $\mathcal Q_0$ from $\mathcal Q$
		\EndIf
		\State $t' \longleftarrow t' + 1$
		\EndWhile
		\State // \textit{Schedule the job and calculate the performance}
		\State $\sigma \longleftarrow \mathcal Q_0$
		\State Remove $\mathcal Q_0$ from $\mathcal Q$
		\State $f_{\lfloor r_{sj}/L\rfloor} \longleftarrow f_{\lfloor r_{sj}/L\rfloor} + \max\left\{ 0, \sigma - r_{sj} \right\}$
		\State Add $\max\{ \sigma, r_{sj} \} + \rho_{sj}$ to $\mathcal Q$
		\EndFor \\
		\textbf{Return:} $f_1,\ldots,f_H$
	\end{algorithmic}
\end{algorithm}

In algorithm \ref{alg:scheduling:main}, $\mathcal Q$ is a vector whose length changes over the time horizon. The entries of $\mathcal Q$ represent the active agents, and the next points in time that each agent becomes available. We initialise $\mathcal Q$ as the zero vector with $y_1'$ entries. We let $\mathcal Q_0$ denote the smallest entry of $\mathcal Q$. In our Python implementations we model $\mathcal Q$ using a heap queue. The simulation iterates through the jobs. In each iteration we advance in time, adding or removing staff as necessary, until one is available to process the current job. We then calculate the delay of that job and update $\mathcal Q$.

We can accelerate the simulation considerably by introducing some specific data structures. First and foremost we can cache previously-computed performance measures. Over the course of the entire Benders decomposition algorithm, we will need the performance measures of identical staffing vectors more than once. By caching the simulation, we can avoid a considerable amount of duplicate computation. 

In a pre-processing phase we do an initial run of the algorithm using the maximum number of staff available in all time periods. Over the course of this initial run, we can save several useful results. For each time period $t \in T$ we can store the index of the first job released in time period $t$ in sample $s$, and the current state of the queue $\mathcal Q$ at the start of that time period. If in a future iteration, $\hat t$ is the first time period with $y_{\hat t}' < M$, then the first $\hat t - 1$ time periods of the simulation are identical to the initialisation. We can save duplicate computation by loading the current state of the queue, and only scheduling jobs which arrive from that time period on. Finally, we can terminate the simulation once $t' > \max\mathcal L$.

In the case of the ACCA problem we also want to impose a penalty of $\lvert T\rvert \times L$ if a passenger is not served by their deadline. Similarly in the NHSS problem, we can impose a penalty for jobs that are not finished by the end of the night shift. These are straightforward changes.

%% file: lbbdpaper.bbl
\begin{thebibliography}{33}
\expandafter\ifx\csname natexlab\endcsname\relax\def\natexlab#1{#1}\fi
\providecommand{\url}[1]{\texttt{#1}}
\providecommand{\href}[2]{#2}
\providecommand{\path}[1]{#1}
\providecommand{\DOIprefix}{doi:}
\providecommand{\ArXivprefix}{arXiv:}
\providecommand{\URLprefix}{URL: }
\providecommand{\Pubmedprefix}{pmid:}
\providecommand{\doi}[1]{\href{http://dx.doi.org/#1}{\path{#1}}}
\providecommand{\Pubmed}[1]{\href{pmid:#1}{\path{#1}}}
\providecommand{\bibinfo}[2]{#2}
\ifx\xfnm\relax \def\xfnm[#1]{\unskip,\space#1}\fi
\bibitem[{Araujo and Repolho(2015)}]{AraujoRepolho2015}
\bibinfo{author}{Araujo, G.}, \bibinfo{author}{Repolho, H.},
  \bibinfo{year}{2015}.
\newblock \bibinfo{title}{Optimizing the airport check-in counter allocation
  problem}.
\newblock \bibinfo{journal}{Journal of Transport Literature}
  \bibinfo{volume}{9}, \bibinfo{pages}{15--19}.
\bibitem[{Beck(2010)}]{Beck2010}
\bibinfo{author}{Beck, J.C.}, \bibinfo{year}{2010}.
\newblock \bibinfo{title}{Checking-up on branch-and-check}, in:
  \bibinfo{booktitle}{Principles and Practice of Constraint Programming -- CP
  2010}, \bibinfo{publisher}{Springer Berlin Heidelberg}.
\bibitem[{Bekker et~al.(2019)Bekker, Moeke and
  Schmidt}]{MoekeBekkerSchmidt2019}
\bibinfo{author}{Bekker, R.}, \bibinfo{author}{Moeke, D.},
  \bibinfo{author}{Schmidt, B.}, \bibinfo{year}{2019}.
\newblock \bibinfo{title}{Keeping pace with the ebbs and flows in daily nursing
  home operations}.
\newblock \bibinfo{journal}{Health care management science}
  \bibinfo{volume}{22}, \bibinfo{pages}{350--363}.
\bibitem[{Benders(1962)}]{Benders1962}
\bibinfo{author}{Benders, J.}, \bibinfo{year}{1962}.
\newblock \bibinfo{title}{{Partitioning procedures for solving mixed-variables
  programming problems}}.
\newblock \bibinfo{journal}{Numerische Mathematik} \bibinfo{volume}{4},
  \bibinfo{pages}{238--252}.
\bibitem[{Bruno et~al.(2019)Bruno, Diglio, Genovese and
  Piccolo}]{BrunoDiglioGenovesePicollo2019}
\bibinfo{author}{Bruno, G.}, \bibinfo{author}{Diglio, A.},
  \bibinfo{author}{Genovese, A.}, \bibinfo{author}{Piccolo, C.},
  \bibinfo{year}{2019}.
\newblock \bibinfo{title}{A decision support system to improve performances of
  airport check-in services}.
\newblock \bibinfo{journal}{Soft Computing} \bibinfo{volume}{23}.
\bibitem[{Bruno and Genovese(2010)}]{BrunoGenovese2010}
\bibinfo{author}{Bruno, G.}, \bibinfo{author}{Genovese, A.},
  \bibinfo{year}{2010}.
\newblock \bibinfo{title}{A mathematical model for the optimization of the
  airport check-in service problem}.
\newblock \bibinfo{journal}{Electron. Notes Discret. Math.}
  \bibinfo{volume}{36}, \bibinfo{pages}{703--710}.
\bibitem[{Codato and Fischetti(2006)}]{CodatoFischetti2006}
\bibinfo{author}{Codato, G.}, \bibinfo{author}{Fischetti, M.},
  \bibinfo{year}{2006}.
\newblock \bibinfo{title}{Combinatorial {Benders'} cuts for mixed-integer
  linear programming}.
\newblock \bibinfo{journal}{Operations Research} \bibinfo{volume}{54},
  \bibinfo{pages}{756--766}.
\bibitem[{Dieleman et~al.(2019)Dieleman, Buitink, den Uijl, Koreman, Passial,
  Couwenberg, Bekker, Moeke and Otsen}]{BekkerEtAl2019}
\bibinfo{author}{Dieleman, N.}, \bibinfo{author}{Buitink, M.},
  \bibinfo{author}{den Uijl, J.}, \bibinfo{author}{Koreman, K.},
  \bibinfo{author}{Passial, R.}, \bibinfo{author}{Couwenberg, M.},
  \bibinfo{author}{Bekker, R.}, \bibinfo{author}{Moeke, D.},
  \bibinfo{author}{Otsen, F.}, \bibinfo{year}{2019}.
\newblock \bibinfo{title}{Demand-driven task-scheduling in a nursing home
  setting: A genetic algorithm approach}.
\newblock \bibinfo{journal}{SSRN Electronic Journal} .
\bibitem[{van Eeden et~al.(2016)van Eeden, Moeke and
  Bekker}]{vanEedenMoekeBekker2016}
\bibinfo{author}{van Eeden, K.}, \bibinfo{author}{Moeke, D.},
  \bibinfo{author}{Bekker, R.}, \bibinfo{year}{2016}.
\newblock \bibinfo{title}{Care on demand in nursing homes: a queueing theoretic
  approach}.
\newblock \bibinfo{journal}{Health care management science}
  \bibinfo{volume}{19}, \bibinfo{pages}{227--240}.
\bibitem[{{Elci} and {Hooker}(2022)}]{OzgunHooker2020}
\bibinfo{author}{{Elci}, O.}, \bibinfo{author}{{Hooker}, J.N.},
  \bibinfo{year}{2022}.
\newblock \bibinfo{title}{{Stochastic Planning and Scheduling with Logic-Based
  Benders Decomposition}}.
\newblock \bibinfo{journal}{INFORMS Journal on Computing} \bibinfo{volume}{0}.
\bibitem[{Fazel-Zarandi et~al.(2012)Fazel-Zarandi, Berman and
  Beck}]{ZarandiBermanBeck2012}
\bibinfo{author}{Fazel-Zarandi, M.}, \bibinfo{author}{Berman, O.},
  \bibinfo{author}{Beck, C.}, \bibinfo{year}{2012}.
\newblock \bibinfo{title}{Solving a stochastic facility location/fleet
  management problem with logic-based benders decomposition}.
\newblock \bibinfo{journal}{Iie Transactions} \bibinfo{volume}{45}.
\bibitem[{Fischetti et~al.(2019)Fischetti, Ljubi\'{c}, Monaci and
  Sinnl}]{Fischetti2019}
\bibinfo{author}{Fischetti, M.}, \bibinfo{author}{Ljubi\'{c}, I.},
  \bibinfo{author}{Monaci, M.}, \bibinfo{author}{Sinnl, M.},
  \bibinfo{year}{2019}.
\newblock \bibinfo{title}{Interdiction games and monotonicity, with application
  to knapsack problems}.
\newblock \bibinfo{journal}{INFORMS Journal on Computing} \bibinfo{volume}{31},
  \bibinfo{pages}{390--410}.
\bibitem[{Geoffrion(1972)}]{Geoffrion1972}
\bibinfo{author}{Geoffrion, A.}, \bibinfo{year}{1972}.
\newblock \bibinfo{title}{Generalized {Benders} decomposition}.
\newblock \bibinfo{journal}{Journal of Optimization Theory and Applications}
  \bibinfo{volume}{10}, \bibinfo{pages}{237--260}.
\bibitem[{Guo et~al.(2021)Guo, Bodur, Aleman and
  Urbach}]{GuoBodurAlemanUrbach2019}
\bibinfo{author}{Guo, C.}, \bibinfo{author}{Bodur, M.},
  \bibinfo{author}{Aleman, D.M.}, \bibinfo{author}{Urbach, D.R.},
  \bibinfo{year}{2021}.
\newblock \bibinfo{title}{Logic-based benders decomposition and binary decision
  diagram based approaches for stochastic distributed operating room
  scheduling}.
\newblock \bibinfo{journal}{INFORMS Journal on Computing} \bibinfo{volume}{0}.
\bibitem[{{Gurobi Optimization, LLC}(2022)}]{Gurobi}
\bibinfo{author}{{Gurobi Optimization, LLC}}, \bibinfo{year}{2022}.
\newblock \bibinfo{title}{{Gurobi Optimizer Reference Manual}}.
\bibitem[{Hooker(2000)}]{Hooker2000}
\bibinfo{author}{Hooker, J.N.}, \bibinfo{year}{2000}.
\newblock \bibinfo{title}{Logic-based Methods for Optimization: Combining
  Optimization and Constraint Satisfaction}.
\newblock \bibinfo{publisher}{Wiley, New York}.
\bibitem[{Hooker(2019)}]{Hooker2019}
\bibinfo{author}{Hooker, J.N.}, \bibinfo{year}{2019}.
\newblock \bibinfo{title}{Logic-Based Benders Decomposition for Large-Scale
  Optimization}. \bibinfo{publisher}{Springer International Publishing}.
\newblock pp. \bibinfo{pages}{1--26}.
\bibitem[{Hooker and Ottosson(2003)}]{HookerOttosson2003}
\bibinfo{author}{Hooker, J.N.}, \bibinfo{author}{Ottosson, G.},
  \bibinfo{year}{2003}.
\newblock \bibinfo{title}{Logic-based benders decomposition}.
\newblock \bibinfo{journal}{Mathematical Programming} \bibinfo{volume}{96},
  \bibinfo{pages}{33--60}.
\bibitem[{Lalita et~al.(2020)Lalita, Manna and Murthy}]{LalitaMannaMurthy2020}
\bibinfo{author}{Lalita, T.R.}, \bibinfo{author}{Manna, D.K.},
  \bibinfo{author}{Murthy, G.S.R.}, \bibinfo{year}{2020}.
\newblock \bibinfo{title}{Mathematical formulations for large scale check-in
  counter allocation problem}.
\newblock \bibinfo{journal}{Journal of Air Transport Management}
  \bibinfo{volume}{85}, \bibinfo{pages}{101796}.
\bibitem[{Laporte and Louveaux(1993)}]{LaporteLouveaux1993}
\bibinfo{author}{Laporte, G.}, \bibinfo{author}{Louveaux, F.V.},
  \bibinfo{year}{1993}.
\newblock \bibinfo{title}{The integer {L-Shaped} method for stochastic integer
  programs with complete recourse}.
\newblock \bibinfo{journal}{Operations Research Letters} \bibinfo{volume}{13},
  \bibinfo{pages}{133--142}.
\bibitem[{Lieder et~al.(2015)Lieder, Moeke, Koole and
  Stolletz}]{LiederMoekeKooleStolletz2015}
\bibinfo{author}{Lieder, A.}, \bibinfo{author}{Moeke, D.},
  \bibinfo{author}{Koole, G.}, \bibinfo{author}{Stolletz, R.},
  \bibinfo{year}{2015}.
\newblock \bibinfo{title}{Task scheduling in long-term care facilities: A
  client-centered approach}.
\newblock \bibinfo{journal}{Operations Research for Health Care}
  \bibinfo{volume}{6}, \bibinfo{pages}{11--17}.
\bibitem[{Lombardi et~al.(2010)Lombardi, Milano, Ruggiero and
  Benini}]{LombardiMilanoRuggieroBenini2010}
\bibinfo{author}{Lombardi, M.}, \bibinfo{author}{Milano, M.},
  \bibinfo{author}{Ruggiero, M.}, \bibinfo{author}{Benini, L.},
  \bibinfo{year}{2010}.
\newblock \bibinfo{title}{Stochastic allocation and scheduling for conditional
  task graphs in multi-processor systems-on-chip}.
\newblock \bibinfo{journal}{Journal of Scheduling} \bibinfo{volume}{13},
  \bibinfo{pages}{315--345}.
\bibitem[{Moeke et~al.(2016)Moeke, van~de Geer, Koole and
  Bekker}]{MoekeGeerKooleBekker2016}
\bibinfo{author}{Moeke, D.}, \bibinfo{author}{van~de Geer, R.},
  \bibinfo{author}{Koole, G.}, \bibinfo{author}{Bekker, R.},
  \bibinfo{year}{2016}.
\newblock \bibinfo{title}{On the performance of small-scale living facilities
  in nursing homes: a simulation approach}.
\newblock \bibinfo{journal}{Operations research for health care}
  \bibinfo{volume}{11}, \bibinfo{pages}{20--34}.
\bibitem[{Moeke et~al.(2014)Moeke, Koole and
  Verkooijen}]{MoekeKoolVerkooijen2014}
\bibinfo{author}{Moeke, D.}, \bibinfo{author}{Koole, G.},
  \bibinfo{author}{Verkooijen, L.}, \bibinfo{year}{2014}.
\newblock \bibinfo{title}{Scale and skill-mix efficiencies in nursing home
  staffing: inside the black box}.
\newblock \bibinfo{journal}{Health Systems} \bibinfo{volume}{3},
  \bibinfo{pages}{18--28}.
\bibitem[{Naderi and Roshanaei(2022)}]{NaderiRoshanaei2021}
\bibinfo{author}{Naderi, B.}, \bibinfo{author}{Roshanaei, V.},
  \bibinfo{year}{2022}.
\newblock \bibinfo{title}{Critical-path-search logic-based benders
  decomposition approaches for flexible job shop scheduling}.
\newblock \bibinfo{journal}{INFORMS Journal on Optimization}
  \bibinfo{volume}{4}, \bibinfo{pages}{1--28}.
\bibitem[{Pearce and Forbes(2018)}]{PearceForbes2018}
\bibinfo{author}{Pearce, R.H.}, \bibinfo{author}{Forbes, M.},
  \bibinfo{year}{2018}.
\newblock \bibinfo{title}{Disaggregated benders decomposition and
  branch-and-cut for solving the budget-constrained dynamic uncapacitated
  facility location and network design problem}.
\newblock \bibinfo{journal}{European Journal of Operational Research}
  \bibinfo{volume}{270}, \bibinfo{pages}{78--88}.
\bibitem[{Rahmaniani et~al.(2017)Rahmaniani, Crainic, Gendreau and
  Rei}]{Rahmaniani2017}
\bibinfo{author}{Rahmaniani, R.}, \bibinfo{author}{Crainic, T.G.},
  \bibinfo{author}{Gendreau, M.}, \bibinfo{author}{Rei, W.},
  \bibinfo{year}{2017}.
\newblock \bibinfo{title}{The benders decomposition algorithm: A literature
  review}.
\newblock \bibinfo{journal}{European Journal of Operational Research}
  \bibinfo{volume}{259}, \bibinfo{pages}{801--817}.
\bibitem[{van Slyke and Wets(1969)}]{SlykeWets1969}
\bibinfo{author}{van Slyke, R.M.}, \bibinfo{author}{Wets, R.},
  \bibinfo{year}{1969}.
\newblock \bibinfo{title}{{L-Shaped} linear programs with applications to
  optimal control and stochastic programming}.
\newblock \bibinfo{journal}{SIAM Journal on Applied Mathematics}
  \bibinfo{volume}{17}, \bibinfo{pages}{638--663}.
\bibitem[{Thorsteinsson(2001)}]{Erlendur2001}
\bibinfo{author}{Thorsteinsson, E.}, \bibinfo{year}{2001}.
\newblock \bibinfo{title}{Branch-and-check: A hybrid framework integrating
  mixed integer programming and constraint logic programming}.
\newblock \bibinfo{journal}{Lecture Notes in Computer Science}
  \bibinfo{volume}{2239}.
\bibitem[{{van Dijk} and {van der Sluis}(2006)}]{DijkSluis2006}
\bibinfo{author}{{van Dijk}, N.M.}, \bibinfo{author}{{van der Sluis}, E.},
  \bibinfo{year}{2006}.
\newblock \bibinfo{title}{Check-in computation and optimization by simulation
  and ip in combination}.
\newblock \bibinfo{journal}{European Journal of Operational Research}
  \bibinfo{volume}{171}, \bibinfo{pages}{1152--1168}.
\bibitem[{Zhang et~al.(2017)Zhang, Matta, Alfieri and
  Pedrielli}]{ZhangMattaAlfieriPedrielli2017}
\bibinfo{author}{Zhang, M.}, \bibinfo{author}{Matta, A.},
  \bibinfo{author}{Alfieri, A.}, \bibinfo{author}{Pedrielli, G.},
  \bibinfo{year}{2017}.
\newblock \bibinfo{title}{A simulation-based benders' cuts generation for the
  joint workstation, workload and buffer allocation problem}, in:
  \bibinfo{booktitle}{2017 13th IEEE Conference on Automation Science and
  Engineering (CASE)}, pp. \bibinfo{pages}{1067--1072}.
\bibitem[{Zhang et~al.(2018)Zhang, Matta, Alfieri and
  Pedrielli}]{ZhangMattaAlfieriPedrielli2018}
\bibinfo{author}{Zhang, M.}, \bibinfo{author}{Matta, A.},
  \bibinfo{author}{Alfieri, A.}, \bibinfo{author}{Pedrielli, G.},
  \bibinfo{year}{2018}.
\newblock \bibinfo{title}{Simulation-based benders cuts: A new cutting approach
  to approximately solve simulation-optimization problems}, in:
  \bibinfo{booktitle}{Proceedings of the 2018 Winter Simulation Conference},
  \bibinfo{publisher}{IEEE Press}. p. \bibinfo{pages}{2225–2236}.
\bibitem[{Zhang et~al.(2019)Zhang, Matta, Alfieri and
  Pedrielli}]{ZhangMattaAlfieriPedrielli2019}
\bibinfo{author}{Zhang, M.}, \bibinfo{author}{Matta, A.},
  \bibinfo{author}{Alfieri, A.}, \bibinfo{author}{Pedrielli, G.},
  \bibinfo{year}{2019}.
\newblock \bibinfo{title}{Feasibility cut generation by simulation: Server
  allocation in serial-parallel manufacturing systems}, in:
  \bibinfo{booktitle}{Proceedings of the 2019 Winter Simulation Conference},
  \bibinfo{publisher}{IEEE Press}. pp. \bibinfo{pages}{3633--3644}.

\end{thebibliography}
